\documentclass[preprint,12pt]{elsarticle}

\usepackage[papersize={210mm,297mm},lmargin=2.5cm,rmargin=2.5cm,top=3cm,bottom=3cm]{geometry}

\usepackage{amsmath,amsfonts,amssymb,amsthm}
\usepackage{latexsym,pdfsync,xcolor,graphicx}
\usepackage[nottoc]{tocbibind}
\usepackage[T1]{fontenc}
\usepackage[utf8]{inputenc}												
\usepackage[english]{babel}			
\usepackage{tikz}
\usepackage{appendix}
\usepackage{hyperref}
\usepackage{enumerate}
\usepackage[shortlabels]{enumitem}
\usepackage{longtable}
\usepackage{float}
\usepackage{multirow,array}
\usepackage{verbatim}
\usepackage{titlesec}
\usepackage{graphicx}
\usepackage{caption}
\usepackage{subcaption}
\usepackage{multicol}
\usepackage{pifont}

\definecolor{verde}{RGB}{67, 155, 88}
\newcommand{\cmark}{\color{verde}{\ding{51}}}%
\newcommand{\xmark}{\color{red}{\ding{55}}}%

\numberwithin{equation}{section}
\newtheorem*{conjecture*}{Conjecture}

\journal{****}

\begin{document}

\begin{frontmatter}

\title{Predator - prey models: a review on some recent adaptations}
\author[a]{\'Erika Diz-Pita}
\address[a]{Departamento de Estat\'istica, An\'alise Matem\'atica e Optimizaci\'on, Universidade
	de Santiago de Compostela, 15782 Santiago de Compostela, Spain }
\ead{erikadiz.pita@usc.es}
\author[a]{M. Victoria Otero-Espinar}
\ead{mvictoria.otero@usc.es}
\begin{abstract}
	
In the last years, predator-prey systems have increased their applications and have given rise to systems which represent more accurately different biological issues that appear in the context of interacting species. Our aim in this paper is to give a state-of-art review of recent predator-prey models which include some interesting characteristics as Allee effect, fear effect, cannibalism and immigration. We compare the qualitative results obtained for each of them, particularly regarding the equilibria, local and global stability, and the existence of limit cycles.
\end{abstract}
\begin{keyword}
Predator-prey \sep Lotka-Volterra \sep Allee effect \sep fear effect \sep cannibalism \sep immigration
\end{keyword}

\end{frontmatter}

\section{Introduction}
Population dynamics is one of the most widely discussed topics within biomathematics. The study of the evolution of different populations has always been of special interest, starting with populations of a single species, but evolving to more realistic models where different species live and interact in the same habitat. Among them we can find models that study competitive relationships, symbiosis, commensalism or predator-prey relationships.

We will focus on predator-prey systems, and our aim is to give a state-of-art review of recent predator-prey models which include Allee efffect, fear effect, cannibalism and immigration. We compare the qualitative results obtained for each of them, particularly regarding the equilibria, local and global stability and the existence of limit cycles.

Probably the most famous predator-prey models are Lotka-Volterra systems, which are polynomial differential equations of degree two, that were initially proposed by Alfred J. Lotka \cite{Lotka} in 1925 and Vito Volterra \cite{Volterra} in 1926, both in the context of competing species. 

At first, Lotka proposed a model for the study of an herbivorous animal species and a plant species, and then he developed the application to the study of the dynamics of a predator–prey system. On his behalf, Volterra considered the same model simultaneously to explain the observations 
about the percentage of predatory fish caught during the years of World War I, which had increased. This fact seemed confusing, because the fishing effort had been reduced during that years, so Volterra was interested in studying that situation.

In order to explain this behaviour, Volterra stated a simple system of ordinary differential equations. He considered $x(t)$ and $y(t)$ the densities of prey and predators, 
and assuming linear behaviours and taking positive constants $a$, $b$, $c$ and $d$, he stated the following model 
\begin{align*}
	\dot{x}&=x(a-by),\\
	\dot{y}&=y(-c+dx).
\end{align*}

Later on, Lotka-Volterra systems were generalized and considered in arbitrary dimension, i.e.
\begin{equation}\label{LVgen}
	\dot{x_i}=x_i  \left( a_{i0} + \sum_{j=1}^{n} a_{ij}x_j \right)  , \;\;\; i=1,...,n.
\end{equation}
Consequently the applications of these systems started to multiply. New applications on population dynamics had been developed,  but also these systems have been used for modelling many natural phenomena, such as chemical reactions \cite{chemical}, plasma physics \cite{plasma} or hydrodynamics \cite{hidrodinamica} just as other problems from social science and
economics \cite{economics}. In view of this, it is clear that population models, in particular competition or predator-prey models, are important both in their original field and in their application to many other problems from other areas. Their study is also interesting from a theoretical point of view. 

Actually, the study of predator-prey models has always been one of the hot spots in biomathematics, so there is a very big number of works on this topic, especially in dimension two, it is, with one predator species and one prey species. For this reason we believe that this review may be of particular interest to those researchers involved in modelling this kind of biological and ecological problems, allowing them to have an overview of the work who is currently being carried out.

Therefore, our objective in this paper is to analyse several adaptations of predator-prey models that have been done in the last years in order to model different real situations that appear in the field of population dynamics. 

We have decided to focus on some special issues to which researchers have devoted special attention in recent years, particularly  the study of the Allee effect 
, the influence of fear effect 
, cannibalism 
and immigration 
. We introduce some papers with influences between them, which allow us to understand the importance of these biological aspects, through the study of several models that are completed and improved in successive works.

The Allee effect appears when a population of prey has a really small density, so it makes difficult for them to reproduce or survive. It was first introduced by Allee in \cite{Allee}. There are several works in the literature that analyse this effect in different population models and conclude that it can have important effects on the system dynamics, including the stabilisation or destabilisation of a system. It also can cause that the solutions of a system take a much longer time to reach a stable equilibrium point.

Another issue to consider is that in certain ecosystems, prey may feel fear of predators and act accordingly making hunting more difficult for the predators.   The theoretical reasonings about the effect of fear behaviors are supported by real experiments. For example, Zanette , White, Allen and Clinchy \cite{songbirds} conducted an experiment on song sparrows during a whole breeding season by using electrical fence to eliminate direct predation of both juvenile and adult song sparrows. No direct killing can happen in the experiment; however, broadcast of vocal cues of known predators in the field was employed to mimic predation risk. Two groups of female song sparrows were tested, among which one group was exposed to predator sounds while the other group was not. The authors found that the group of song sparrows exposed to predator sounds produced 40\% less offspring than the other group. They believe this is because fewer eggs were laid, fewer eggs were successfully hatched, and fewer nestlings survived eventually. Also because there were behavioral changes including less time of adult song sparrows on brood and less feeding to nestlings during breeding period. Similar experiments on other birds and vertebrate species reported the same conclusions  \cite{creel,sheriff,wirsing}, it is,  even though there is no direct killing between predators and prey, the presence of predators cause a reduction in prey population due to anti-predator behaviors.

We analyze different models with the presence of fear effect, combined with different functional responses, with both omnivorous and specialist predators,  with hunting cooperation, with Allee effect and with prey refuge.

As we mentioned, we will study ecosystems in which predators feed only on one prey species, which we will call specialists, and others in which predators are omnivorous, it is, they will be able to feed themselves from other resources. A special case of this is cannibalism, it is, the act of killing and at least partial consumption of conspecifics. Cannibalism actively occurs in more than 1300 species in the nature \cite{1300can}, and it has been mathematically modelled for some ecosystems, as for example the European regional seas ecosystem model (ERSEM) for the North Sea \cite{ERSEM}.

In the literature, one of the first contributions to the study of this topic was the work of Kohlmeir \cite{Kohlmeier}, who considered cannibalism  in the predator in the classic Rosenzweig-McArthur model. In general, cannibalism has been considered primarily in predators, and the results agree that it has stabilising effects, and can cause the survival of a species that would otherwise be driven to extinction. We have selected two works that illustrate well the different effects of cannibalism when we add it to already well-known models, as the Lotka-Volterra and the Holling-Tanner model. In addition, the inclusion of cannibalism is considered not only in the predator species as usual, but also in the prey species. 

Just as cannibalism or omnivorism can be adaptation strategies to the lack or shortage of food, other types of strategies or behaviours can be induced by this lack of food or by the hostility of the habitat, such as migration. 
In fact, most predator-prey systems in the wild are not isolated, so it is important to consider the effects of the presence of some number of immigrants. We will focus on studying how the inclusion of immigration changes the dynamics of the classic Rosenzweig-MacArthur model and the Lotka-Volterra model.

Although we have chosen the previous topics as the thread of this review, in many of the considered works these ones appear combined with other biological characteristics such as prey refugee, hunting cooperation or with different types of functional reponses. Down below we make some considerations and give some references about this topics.

Firstly, we recall that in ecology a functional response is the intake rate of a predator as a function of food density. The most usual classification of functional responses is the one given by C. S. Holling, in which three types of responses are considered. The type I functional response 
assumes that the consuming and hunting rates are linear up to a maximum where they become constant. This linear increase assumes that the time needed by the predator to process food is negligible and that consuming food does not interfere with searching food. This is the functional response used in the Lotka-Volterra model. 

The type II functional response is given by $f(x)=a x/(1+a h x)$, where $x$ denotes the food or prey density, $a$ is the rate at which the predator encounters food items per unit of food density, and $h$ is the average time spent on processing a food item. This functional response assumes that the consumption rate is decelerated with increased population as a consequence of a limited capacity in searching and processing food. This is the functional response considered in the Rosenzweig-MacArthur model, which is a classical model with logistic growth \cite{RosMac}.

In type III functional response the consumption rate is more than linear at low levels of resource. This is suitable for instance in a population of predators that has to learn how to hunt efficiently. This functional response is given by $f(x)=nx^p/(a^p+x^p)$, where $a$ and $n$ are positive and $p$ is an integer with $p>1$.

These three types of functional responses are widely used in population dynamics models and specially in predator-prey models, not only on ODE models but also on discrete models, diffusion models, stochastic or fractional order models. Some interesting recent works that study this kind of functional responses are \cite{fr1,fr2,fr3,fr4,fr5,fr6}. 

There are other types of functional responses, for example the one proposed by Taylor, also known sometimes as Holling type IV functional response, which is given by the function
\begin{equation}
	h(x)=\frac{q x}{x^2+bx+a},
\end{equation}
where $b\geq0$. Other models consider a functional response of ratio-dependent type, in which the consider function is $h(x)=a/x$, with $a\in \mathbb{R}$. 

Regarding the existence of prey refuge, it has also special interest in the study of predator-prey populations, and many scholars have made great achievements in this aspect. There are interesting works on this subject apart from the ones we include in this review, i.e. those related to our main topics. As a suggestion for interested readers we propose the following references \cite{prey1,prey2,prey3,prey4,Kar}.

Finally we would like to note that the third characteristic that we have mentioned, hunting cooperation on predators, has been investigated by many authors in mathematical modelling, but most of the time independently of the effect of fear on prey \cite{HC1,Berec,HC3,HC4,HC5}. In this review we focus on the combination of both hunting cooperation and fear effect, which we believe is particularly important due to the ecological evidence. For example, wolves cooperate during hunting and when they are present, elks use anti-predator stategies and avoid areas frequented by the wolves \cite{wolves}. Elks respond to the presence of wolves by altering foraging, vigilance, habitat selection, patterns of aggregationg and sensitivity to environment \cite{elks1,elks2,elks3,elks4}. In the same way, while lionesses show cooperative hunting behaviours \cite{lioness}, zebras are affected by fear of predation risk, reaching areas where the encounter with lionesses is less frequent \cite{zebras}.

In the following, we will focus on models with Allee effect in Section \ref{sec:Allee}, with fear effect in Section \ref{sec:fear}, with cannibalism in Section \ref{sec:cannibalism} and with immigration in Section \ref{sec:immigration}. We compare the proposed models whenever it is possible, and we include figures that allow us to compare different models and to see how the inclusion of different biological characteristics affects the dynamics. All these figures have been made with the free software R.  In Table \ref{tab:resumen} we summarise the most important characteristics of the models studied.

\begin{table}[H]
	\begin{center}
					\resizebox{\textwidth}{!}{
	\begin{tabular}{|ccccccc|}
			\hline
			\textbf{Model} & \textbf{Allee effect} & \textbf{Fear} & \textbf{Migration}  &  \textbf{Cannibalism} & \textbf{Funct. resp.} & \textbf{Observations}   \\
	\hline
		\hline
	\eqref{sisLV} &  \xmark & \xmark & \xmark & \xmark & Type I  &  \\	
			\hline
			\eqref{LVAllee_prey}&  \textbf{\textcolor{verde}{On prey}} & \xmark & \xmark  & \xmark & Type I  &  \\	
			\hline
			\eqref{AlleeD} &  \textbf{\textcolor{verde}{On predator}} & \xmark & \xmark & \xmark  &Type I  &  \\	
			\hline
			\eqref{LV_AlleeDP} &  \textbf{\textcolor{verde}{On predator}} & \xmark & \xmark & \xmark &Type I   &  \\
			\hline
			\eqref{Alleefg} &  \textbf{\textcolor{verde}{On prey}} & \xmark & \xmark & \xmark &  Types I, II and III   &  \\		
			\hline
			\eqref{Olivares2010} &  \textbf{\textcolor{verde}{On prey}} & \xmark & \xmark & \xmark & Type I   &  \\	
			\hline
			\eqref{Alleedoble} &  \textbf{\textcolor{verde}{Double on prey}} & \xmark & \xmark & \xmark  & Type I  &  \\	
			\hline
			\eqref{Alleesimple} &  \textbf{\textcolor{verde}{Unique on prey}} & \xmark & \xmark & \xmark  &Type  I  &  \\	
			\hline
			\eqref{15}&  \textbf{\textcolor{verde}{On predator}} & \xmark & \xmark & \xmark & Type II   & Omnivorous predator \\	
			\hline
			\eqref{siswang} &\xmark & \cmark & \xmark & \xmark &Types I and II   &  \\	
			\hline
			\eqref{25} &  \xmark & \cmark & \xmark & \xmark & Type I   &  \\	
			\hline
			\eqref{fearcaza} &  \xmark& \cmark  & \xmark & \xmark  & Type II  & Hunting cooperation  \\	
			\hline
			\eqref{27} &  \textbf{\textcolor{verde}{On prey}} & \cmark  & \xmark & \xmark  & Type I  &  \\	
			\hline
			\eqref{aditivo} &  \textbf{\textcolor{verde}{On prey}} & \cmark & \xmark & \xmark  & Type I  &  \\	
			\hline
			\eqref{m-pr} &  \xmark & \cmark  & \xmark & \xmark & Type II   & Prey refugee \\	
			\hline
			\eqref{33} &  \textbf{\textcolor{verde}{On prey}} &\cmark & \xmark & \xmark & Type I  & Prey refugee  \\	
			\hline
			\eqref{34} &  \xmark& \xmark & \xmark &  \textcolor{verde}{On predator}  & Type I  &  \\	
			\hline
			\eqref{HT} & \xmark & \xmark & \xmark &  \textcolor{verde}{On predator}  & Ratio-dependent  &  \\	
			\hline
			\eqref{HTpredator} &  \xmark & \xmark & \xmark & \textcolor{verde}{On predator} & Ratio-dependent    &  \\	
			\hline
			\eqref{HTprey} &  \xmark& \xmark & \xmark & \textcolor{verde}{On prey} & Ratio-dependent  &  \\	
			\hline
			\eqref{inm1} &  \xmark& \xmark & \cmark  & \xmark &Type II   &  \\	
			\hline
			\eqref{40} & \xmark & \xmark &\cmark  & \xmark & Type I   &  \\	
			\hline
			\eqref{41} & \xmark & \xmark &\cmark  & \xmark   & Types I, II and III &  \\	
			\hline
	\end{tabular} }
	\caption{Summary of the characteristics of the revised predator-prey models.} \label{tab:resumen}
		\end{center}
\end{table}

\section{Study of Allee effect }\label{sec:Allee}
The Allee effect is the biological phenomenon of correlation between the population size or density and the growth rate. Roughly speaking, it occurs when the population of a species has a really small density, so it makes difficult for them to reproduce or survive. 
There are two types of Allee effect depending on the nature of density dependence at low densities, the \textit{strong Allee effect} and the \textit{weak Allee effect}. The distinction between these two types is based on whether or not the population exhibits a critical population size or density. A population with a weak Allee effect has a reduced per capita growth rate at lower population density or size, but even at this low population size or density, the population will always exhibit a positive per capita growth rate. On the other hand, a population 
with a strong Allee effect has a critical population size or density under which the population growth rate becomes negative. Then when the population density or size is below this threshold, the population will be drive to extinction.

None of this types of Allee effect was considered in the first classical models as the Lotka-Volterra model.  This phenomenon was first introduced by Allee in \cite{Allee}.  \color{black} There are several works in the literature that analyse this effect in different population models \cite{Alleestrong,AlleeNOLV8,Olivares2010} and conclude that it can have important effects on the system dynamics. 
These effects include stabilising or destabilising a system by changing the local stability of some singularities, as can be seen in \cite{AlleeNOLV1,AlleeNOLV2,AlleeNOLV3,AlleeNOLV4,AlleeNOLV5} or cause that the solutions of the system take a much longer time to reach a stable equilibrium point  \cite{AlleeNOLV3,AlleeNOLV6}.

Mathematically, the Allee effect is usually represented modifying the growth function. The most common modification is to use a multiplicative factor, in which case the equation for a single species is given by 
\begin{equation}
	\dot{x}=r\left( 1-\frac{x}{K}\right) (x-m)x,
\end{equation}
where $r$ is the intrinsec growth rate and $K$ the environmental carrying capacity. It is said that the Allee effect is \textit{strong} if $m>0$ and \textit{weak} if $m\leq0$. For the strong Allee effect, the per capita growth rate is negative for $x\in(0,m)$, which implies the extinction of the population.

In this section, we will focus first on some Lotka-Volterra systems subject to an Allee effect and then we will study the inclusion of this effect in more general population models.

\vspace{0.5cm}

First of all, in \cite{Allee1} the author begins by studying the following Lotka-Volterra system
\begin{equation}\label{sisLV}
	\begin{split}
		\dot{x}&= r x (1-x)- a x y , \\
		\dot{y}&= a y (x-y),
	\end{split}
\end{equation}
where the variables $x$ and $y$ represent the number of individuals of the prey population, and the number of individuals of the predator population, respectively.
We will maintain this notation throughout all the review.

The equilibrium points for this system are $(0,0)$, $(1,0)$ and $(x^*,y^*)=(r/(a+r), r/(a+r))$.  The most interesting singularity from a biological point of view is the one in the positive quadrant, which represents the coexistence of both species. Based on the eigenvalues of the jacobian matrix, the author concludes that the equilibrium on the axes are unstable and the positive equilibrium is locally asymptotically stable. The interest of this work lies in the comparison with the case in which Allee effect is taken into acount. For that, the following system is considered:
\begin{equation}\label{LVAllee_prey}
	\begin{split}
		\dot{x}&= r \alpha(x) x (1-x)- a x y , \\
		\dot{y}&= a y (x-y),
	\end{split}
\end{equation}
where $\alpha(x)=x/(\beta+x)$ represents the Allee effect and $\beta>0$ is called the ``Allee effect constant''. We understand that as $\beta$ increases, the existing Allee effect on the population is stronger, and the population growth of that species decreases.

Based on biological facts the function $\alpha(x)$ must verify some conditions:  $\alpha'(x)$ must be positive for any positive $x$ because the Allee effect decreases as density increases and $\lim_{x\to \infty}\alpha(x)=1$ because the Allee effect vanishes at high densities. The authors also state that $\alpha(0)$ must be zero because reproduction is not possible without any individuals, but as the net growth rate is $x\alpha(x)$, the condition $\alpha(0)=0$ is not necessary as $x\alpha(x)$ is always zero if $x=0$, regardless of the value of $\alpha(0)$.

\color{black}

The proposed function $\alpha(x)=x/(\beta+x)$ verifies this set of conditions, but maybe it would be interesting to  analyze different functions representing the Allee effect.

For system \eqref{LVAllee_prey} there exist again three equilibria. Two of them are the same as in the system with no Allee effect, $(0,0)$ and $(1,0)$, and the third one acquires a new expression $(x_*,y_*)=((r-a\beta)/(a+r),(r-a\beta)/(a+r))$. It is important to recall that now, the equilibrium $(x_*,y_*)$ has only biological meaning if $r-a \beta$ is positive.

The stability of the equilibria is studied by means of the eigenvalues. The author says that the equilibria $(0,0)$ and $(1,0)$ are unstable, but the proof for the origin is not given. Since it is a linearly zero equilibrium, it is, the linear part of the system \eqref{LVAllee_prey} in the origin is the identically-zero matrix, it would be convenient to give  a more detailed analysis. 
The equilibrium $(x_*,y_*)$ is locally asymptotically stable.

Numerical simulations are also carried out in the article. 
The outstanding conclusions of \cite{Allee1}, obtained from both analytical and numerical considerations, are the following: The Allee effect may cause that there is no equilibrium in the positive quadrant, while in the classical model this singularity always exists, and furthermore it is locally asymptotically stable. In the case that the Allee effect does not affect the existence of equilibria, it also has consequences on the dynamics of the model; on the one hand, the system takes much longer time to reach the state of equilibrium, and on the other hand, the values of population density of both prey and predators in that point are reduced with respect to the classical model.

\vspace{0.5cm}
In \cite{Allee2} the study of the case in which the Allee effect affects predator population rather than prey is carried out, based on the reason that higher levels in the food chain are more likely to become extinct. The proposed system is analogous to the one in \cite{Allee1}, considering the same function $\alpha$, but in this case as a function of the density of predators:
\begin{equation}\label{AlleeD}
	\begin{split}
		\dot{x}&=r x (1-x)- a x y, \\
		\dot{y}&= a \alpha(y)y (x-y).
	\end{split}
\end{equation}

In this case, for the study of the stability, different results are combined: the study of the trace and determinant of the jacobian matrix of the system, the Bendixon-Dulac criterion and the property of persistence. 

This notion of persistence or permanence appears frequently in ecology as it has special importance on biological systems because it is a condition ensuring the long-term survival of all species. Considering a general system 
\begin{equation*}
	\dot{x_i}=x_i f_i(x), \;\; i=1,...,n,
\end{equation*}
on the positive region $R^n_+$, where $x=(x_1,...,x_n)$ and conditions ensuring the global existence and uniqueness of solutions in forward time hold, it is called \textit{permanent} or \textit{persistent} if there exist $m,M\in (0,\infty)$ such that, given any $x\in \text{int} R^n_+$, there is a $t_x$ such that 
\begin{equation*}
	m\leq x_i(t)\leq M, \;\; i=1,...,n, \;\; t\geq t_x.
\end{equation*}
From a biological point of view, it is reasonable to expect that a permanent system will have a coexistence equilibrium in $\text{int} R_+^n$, as can be seen in \cite{permanence}.

Because of the importance of this property, the first objective addressed by the authors in \cite{Allee2} is the proof of the persistent property of system \eqref{AlleeD}, because they rely on that to prove the local stability. They conclude that the system is persistent under condition $r>a$.

The persistence property makes it possible to ensure that the equilibria on the axes, $(0,0)$ and $(1,0)$ are unstable, as no positive solution can approach to these equilibria. For the equilibrium point in the positive quadrant $(x^*,y^*)=(r(1-x)/a,r(1-x)/a)$, they conclude that it is locally asymptotically stable, since the real part of the eigenvalues is negative. 
Locally, the stability of the equilibria is the same as the obtained in \cite{Allee1} for the system with Allee effect on the prey species.

Here, beyond the local study of equilibria, a simple argument is used to prove that the asymptotic stability of the point $(x^*,y^*)$ is global. They prove that there are no limit cycles, which is easily deduced applying the Bendixon-Dulac criterion, and that is enough
because we know $(x^*,y^*)$ is the only stable equilibrium.
This global result is not proved for the model in \cite{Allee1}.

\begin{figure}[h]
	\centering
	\includegraphics[width=8cm]{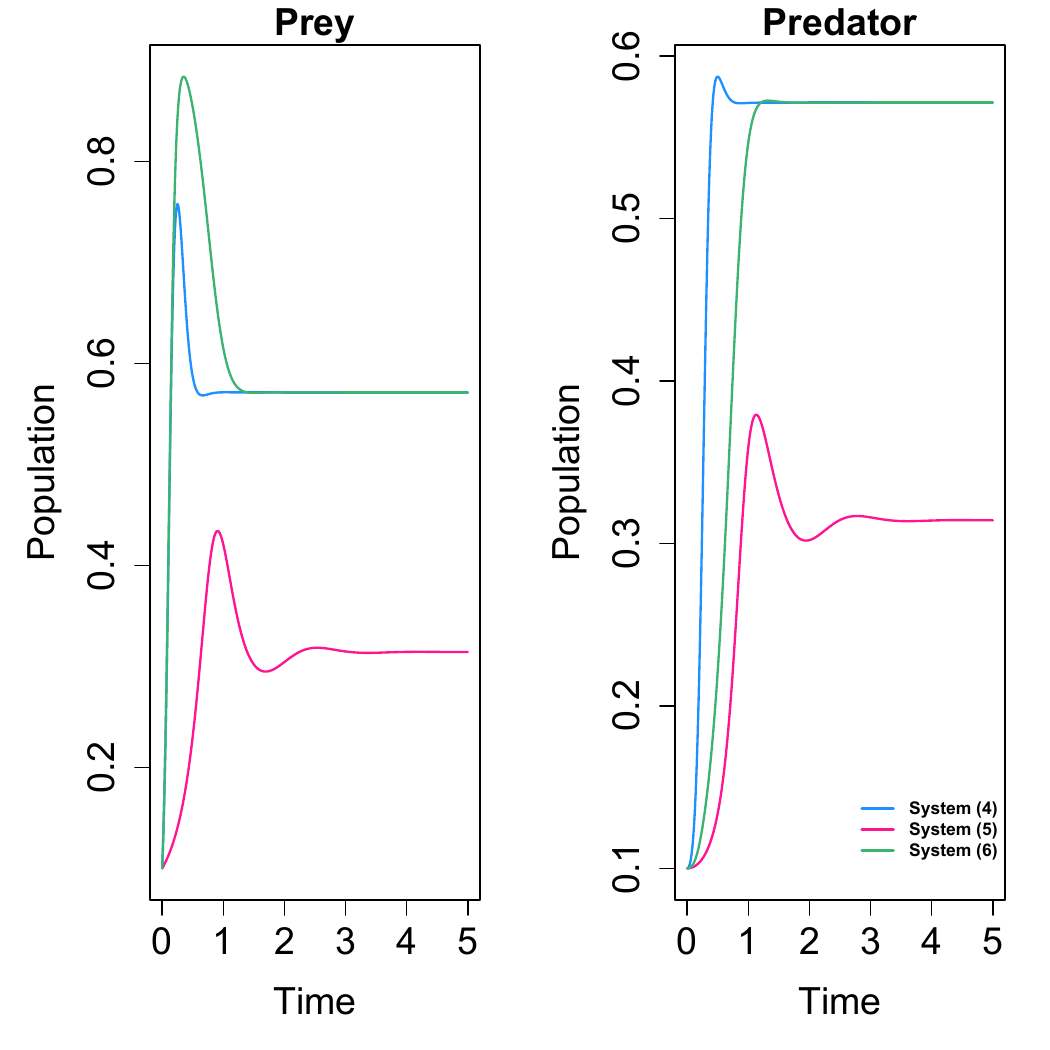}
	\caption{Time–population density of prey and predator for systems \eqref{sisLV}, \eqref{LVAllee_prey} and \eqref{AlleeD}. The initial conditions are $x(0) = 0.1$ and $y(0) = 0.1$, and the values for the parameters are $r = 20$, $a = 15$ and $\beta= 0.6$.}
	\label{fig:Allee1}
\end{figure}

\begin{figure}[h]
	\centering
	\includegraphics[width=8cm]{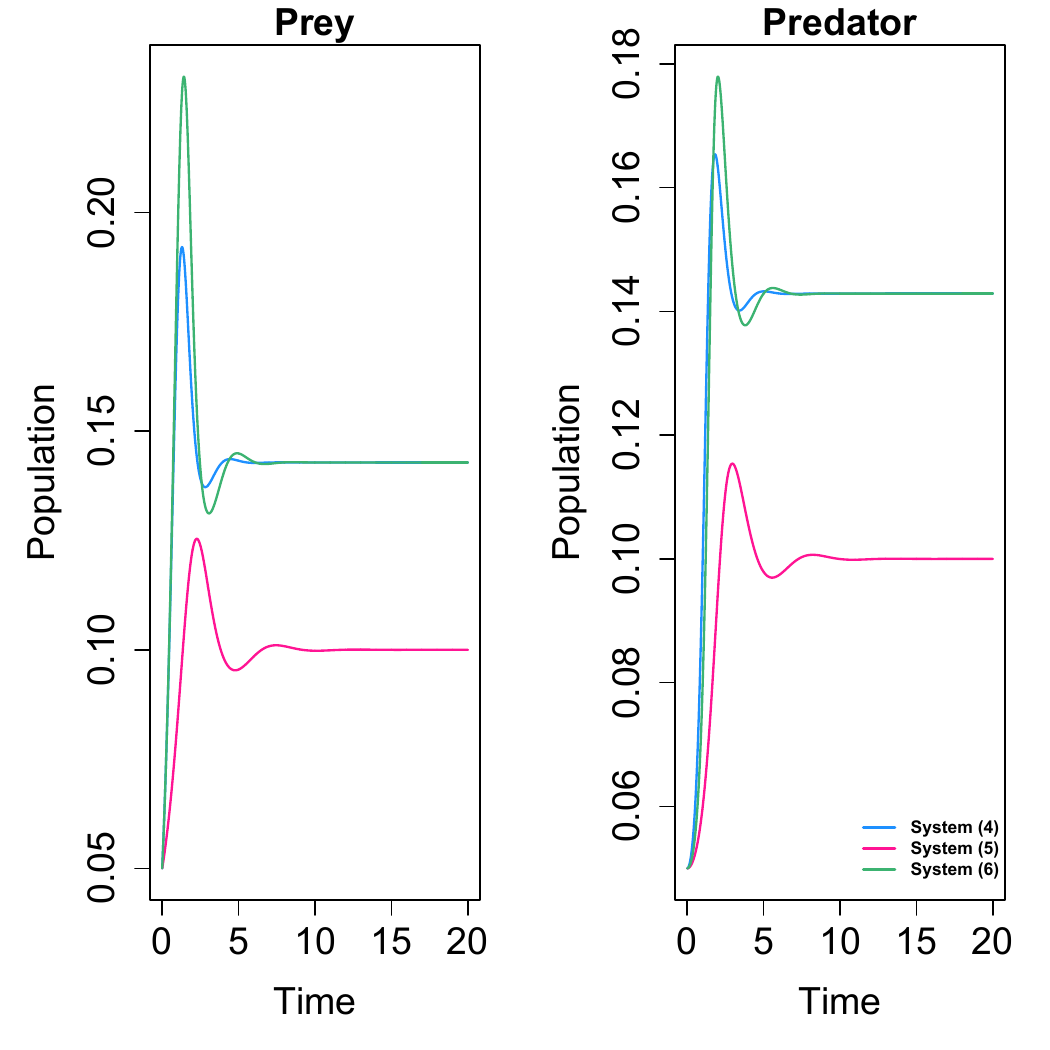}
	\caption{Behaviour of population density of prey and predator for systems \eqref{sisLV}, \eqref{LVAllee_prey} and \eqref{AlleeD}. The initial conditions are $x(0) = 0.005$ and $y(0) = 0.01$, and the values for the parameters are $r = 2.5$, $a = 15$ and $\beta= 0.05$.}
	\label{fig:Allee2}
\end{figure}

After the stability analysis, numerical simulations are carried out in \cite{Allee2}. They show that if $r>a$, the Allee effect on the predator species has no influence on the existence and stability of the positive equilibrium and on the final density of the predator and prey species. This result differs from the obtained in \cite{Allee1}. This is showed in Figure \ref{fig:Allee1}. A similar behaviour is also obtained in some cases with $r<a$, as for example in the case given in Figure 	\ref{fig:Allee2}.

Numerical simulations also show that the system with Allee effect on the predators takes much longer time to reach its positive equilibrium when $\beta$ increases its value, which coincides with the results in the case of Allee effect on the prey species.


\vspace{0.5cm}
In \cite{Allee3} the previous model is combined with the choice of a  more realistic model for the growth of  prey, instead of the logistic one used previously. The proposed model is
\begin{equation}\label{LV_AlleeDP}
	\begin{split}
		\dot{x}&= x \left( \frac{a_1}{b_1+c_1 x}-d_1-e_1 x \right)-axy,\\[0.2cm]
		\dot{y}&= \frac{ay}{\beta+y}y (x-y).
	\end{split}
\end{equation}

Particularly, they take the Beverton-Holt function \cite{Beverton}  as the birth rate.   In the same way that in \cite{Allee2}, the authors prove that the system is permanent. In this case the condition under which the system \eqref{LV_AlleeDP} is permanent is: $a_1/b_1>d_1+a u^*$, where \color{black}
\begin{equation}
	u^*=\frac{-(e_1b_1+d_1c_1)+\sqrt{(e_1b_1+d_1c_1)^2-4e_1c_1(d_1b_1-a_1)}}{2e_1c_1}.
\end{equation} 
As the authors say, this condition can be changed into other which is simpler but more restrictive
\begin{equation}
	\frac{a_1}{b_1}>d_1+ a \frac{(a_1/b_1)-d_1}{e_1}.
\end{equation}

Analogous to the previous system, there exist three equilibria: $(0,0)$, $(u^*,0)$ and $(x^*,y^*)$, and the conclusions about their stability are similar: if the condition that makes the system permanent holds, then the equilibrium points on the axes are unstable, and the positive equilibrium, which appears if $a_1>b_1d_1$, is always locally asymptotically stable, and globally asymptotically stable if the condition for the permanence holds.

Here a extinction property is proved: if $a_1/b_1 <d_1$ then
\begin{equation*}
	\lim_{t \to \infty} x(t)=0 \;\text{and} \; \lim_{t\to \infty} y(t)=0.
\end{equation*}
This kind of property was not proved in the systems studied in \cite{Allee1} and \cite{Allee2}.

At last,  the numerical simulations in this case are only oriented to give some examples for the validation of the conclusions about the stability, but the influence of the Allee effect is not analysed. 

\vspace{0.5cm}

We will now consider the inclusion of Allee effect in more general predator-prey models. In \cite{Alleestrong} the authors study the global dynamics of a model with increasing and bounded functional response which includes a strong Allee effect. They consider a system 
\begin{equation}\label{Alleefg}
	\begin{split}
		\dot{x}&=g(x)(f(x)-y),\\
		\dot{y}&=y (g(x)-d),
	\end{split}
\end{equation}
where the functions $f$ and $g$ verify the following hypothesis:
\begin{itemize}
	\item $f\in \mathcal{C}^1([0,\infty))$, $f(b)=f(1)=0$, with $0<b<1$, $f(x)$ is positive for $b<x<1$ and negative in other case. 
	There exists $\overline{\lambda}\in (b,1)$ such that $f'(x)>0$ in $[b,\overline{\lambda})$ and $f'(x)<0$ in $(\overline{\lambda},1]$.
	\vspace{0.15cm}
	
	\item $g\in \mathcal{C}^1([0,\infty))$, $g(0)=0$, $g(x)>0$ and $g'(x)>0$ if $x>0$. There exists $\lambda>0$ such that $g(\lambda)=d$.
\end{itemize}
The function $g$ is the functional response of predator and $g(x)f(x)$ is the net growth rate of the prey. The conditions about $g$ include Holling type I, II and III functional responses, but no type IV.

The hypothesis on $f$ reflect that in the absence of the predator, the prey has a strong Allee effect growth. 
The equation \eqref{Alleefg} is dimensionless and the conversion efficiency of prey biomass into predator biomass is $1$. The parameter $b$ is the survival threshold or sparsity constant of the prey.
The parameter $d$ is the mortality rate of predator but also $\lambda$, which is the stationary prey population density coexisting with predator, can be thought as a measure of the predator mortality as it increases with $d$. As pointed out by Bazykin \cite{Bazykin}, it is natural to regard $\lambda$ as a measure of how well the predator is adapted to the prey. In \cite{Alleestrong} the parameter $\lambda$ plays an important role as a bifurcation parameter.

A complete theoretical analysis of the model is carried out, including an analysis of the phase portrait and Hopf bifurcation. The system is transformed into other of Liénard type. The existence and uniqueness of limit cycles and heteroclinical orbits are studied, which is important for some ecological phenomena. The first quadrant is shown to be invariant.

There exist four possible equilibria: $(0,0)$, $(1,0)$, $(b,0)$ and $\left( \lambda,f(\lambda)\right) $. This 
latest equilibrium only exists if $b<\lambda<1$, and in any other case the system has three boundary equilibria, there is no coexistence and no limit cycles. In that case the origin is globally asymptotically stable if $0<\lambda\leq b$, and if $\lambda\geq 1$ the regions of attraction of the origin and $(1,0)$ divide the first quadrant.

If $b<\lambda<1$ then the origin is a stable node, $(b,0)$ is a saddle with stable manifold on the $x$-axis and unstable manifold on the prey nullcline $y=f(x)$. The equilibrium $(1,0)$ is a saddle with stable manifold on the $x$-axis and unstable manifold on the region $y>f(x)$. 

The authors also analyse the existence and multiplicity of limit cycles, and for some conditions the uniqueness of limit cycle is proved.

The system has a Hopf bifurcation when $\lambda=\overline{\lambda}$ in $(\overline{\lambda}, f(\overline{\lambda}))$. The Hopf bifurcation can be supercritical and 
backwards (respectively, subcritical and forward) if the first Lyapunov coefficient is negative (respectively positive). Hopf bifurcation is backwards if there exists a periodic orbit of small amplitude for each $\lambda \in (\overline{\lambda}-\varepsilon,\overline{\lambda})$ (similarly with forward), it is supercritical if the periodic solutions which bifurcate are orbitally asymptotically stable (respectively, subcritical if periodic solutions which bifurcate are unstable).

The key bifurcation parameter $\lambda$ is the location of the predator nullcline $x = \lambda$. When $\lambda$ is below the prey extinction threshold $b$, the prey will become extinct and consequently the predator will also starve to death. The authors show that there exists  a threshold $\lambda^*<\overline{\lambda}$  \color{black} such that if $\lambda<\lambda^*$, then the predator invasion leads to the extinction of both species. This already studied phenomenon is known as overexploitation \cite{Rosenz}, and mathematically it means the global stability of the equilibrium $(0,0)$.
The threshold $\lambda^*$ is at the unique point where a point-to-point heteroclinic orbit loop exists.

When $\lambda$ moves across $\lambda=\lambda^*$, a limit cycle appears. It can be established if initially the prey is over the  threshold $b$, and the predator is under the overexploit threshold to make a successful invasion.
There is another change when it moves past $\overline{\lambda}$, as it is the value at which the Hopf bifurcation takes place.  \color{black} Then, the alternative stable state of limit cycle switches to a coexistence equilibrium point, and the oscillation dies down to the stable equilibrium. Finally when $\lambda>1$, the alternative stable state becomes the equilibrium in which there are only prey and the predator species is not able to establish itself and will become extinct. For these cases, higher initial predator density will always lead to the overexploitation.

The authors apply these results to special cases determined by different functional responses $g(u)$, and they study concrete models, setting the expressions of both functions $f$ and $g$. In particular, they apply their results to a cubic model with Holling type II funcional response, a cubic model with linear functional response and 
to the Boukal-Sabelis-Berec model \cite{Sabelis}.

\vspace{0.5cm}

After analysing this more general model, we study another one that can not be written in the form \eqref{Alleefg}.
In \cite{Olivares2010} the studied model is based on the Leslie-Gower model, i. e., the predator growth is given by a logistic-type equation and the environmental carrying capacity of predators is a function on the prey population, $f(x)=nx$. The Leslie-Gower model has only one positive equilibrium which is asymptotically stable. Adding the prey population affected by Allee effect and the lineal functional response, the authors propose the following system:
\begin{equation}\label{Olivares2010}
	\begin{split}
		\dot{x}&=\left(r\left( 1-\frac{x}{K}\right)(x-m)-qy  \right) x,\\[0.1cm]
		\dot{y}&=s\left( 1-\frac{y}{nx}\right)y, 
	\end{split}
\end{equation}
defined in $(0,\infty]\times[0,\infty]$.
The parameter $m$ is the Allee threshold or minimum of viable population, they consider $m\in[0,K)$, resulting in a weak Allee effect when $m=0$ and in a strong Allee effect in any other case. 

The system \eqref{Olivares2010} has always the boundary equilibria $(K,0)$ and $(m,0)$. The positive equilibria $(x^*,nx^*)$ such that $x^*$ is a solution of $rx^2-(Kr+mr-Knq)x+Kmr=0$ there exist if and only if $Kr-mr-Knq>0$ and $\Lambda=(Kr+mr-Knq)^2-4Kmr^2\geq0$. 

By changing the variables by $x=Ku$, $y=nKv$ and rescaling the time $\tau=trK/(u+c/(nK))$, the authors obtain the following Kolmogorov system, which is topologically equivalent and a continuous extension of system \eqref{Olivares2010}, defined at $[0,\infty)\times [0,\infty)$, with $0\leq M<1$:
\begin{equation}\label{Oli10cv}
	\begin{split}
		\frac{du}{d\tau}&=\left( (1-u)(u-M)-Qv\right) u^2,\\
		\frac{dv}{d\tau}&=S(u-v)v.
	\end{split}
\end{equation}
Here $M=m/K$, $Q=qn/r$, $C=c/nK$ and $S=s/rK$.
The equilibrium points of system \eqref{Oli10cv} are $(0,0)$, $(M,0)$, $(1,0)$ and $(u^*,u^*)$.  For this system the authors obtain that for certain values of the parameters there are two positive equilibria $(u_1^*,u_1^*)$ and $(u_2^*,u_2^*)$, only one positive equilibrium or none. 

In the case with strong Allee effect, i. e. with $M>0$, the equilibrium $(M,0)$ is a hyperbolic repellor point and $(1,0)$ is a hyperbolic saddle point. The point $(0,0)$  has a hyperbolic and a parabolic sector, so there exists a separatrix curve in the phase plane that divides the behavior of the trajectories. This implies that the model is highly sensitive to the initial conditions.
Then, this point is an attractor point for certain trajectories and a saddle point for others. In the proof of this result the blow-up technique is applied, which can be found in \cite{Libro,blowup}.

When the system has two positive equilibria, $(u_1^*,u_1^*)$ is a saddle point and $(u_2^*,u_2^*)$ can be a stable point or an unstable equilibrium surrounded by a stable limit cycle. In the transition from one case to another, the equilibrium $(u_2^*,u_2^*)$ is at least a two order weak focus and a stable limit cycle is generated by a Hopf bifurcation.

There exist values of the parameters for which there is a homoclinic curve determined by the stable and unstable manifold of the point $(u_1^*,u_1^*)$ and there is an unstable limit cycle that bifurcates from the homoclinic orbit surrounding point $(u_2^*,u_2^*)$.

For some values of the parameters the equilibria $(u_1^*,u_1^*)$ and $(u_2^*,u_2^*)$ collapse and give rise to a unique equilibrium point in the interior of the positive quadrant. This point can be a non-hyperbolic attractor node, a non-hyperbolic repellor node or a cusp, in which chase there exists a unique trajectory which attains the equilibrium. When this point is a repellor, the origin is a global attractor.

In the case with weak Allee effect, i. e. with $M=0$,  the equilibrium points $(0,0)$ and $(1,0)$ are a non-hyperbolic saddle and a saddle, respectively. Moreover there exists a unique positive equilibrium if and only if $Q<1$. This equilibrium can be a local attractor or a repellor, and in the transition from one case to another, the equilibrium is a weak focus of order one. The equilibrium can be also a repellor focus surrounded by a limit cycle or a repellor node (in which case the origin is globally asymptotically stable).
This work shows that the Allee effect may destabilize a predator–prey model, since the equilibrium points of the system could be changed from stable to unstable.

\vspace{0.5cm}
In the previous works a simple Allee effect was considered, but mathematically it is possible to represent the Allee effect with more than one modification on the system simultaneously. In \cite{Olivares2011} two models are compared, one with a unique Allee effect and other with a double Allee effect on prey. The model with double Allee effect is given by
\begin{equation}\label{Alleedoble}
	\begin{split}
		\dot{x}&=\left( \frac{r}{x+n}\left(1-\frac{x}{K} \right)(x-m)-qy \right)x,\\
		\dot{y}&= (px-c)y,
	\end{split}
\end{equation}
in $[0,\infty]\times[0,\infty]$, where $n<K$ for biological reasons and $m\in(-K,K)$. 
The Allee effect is expressed once in the factor $\alpha(x)=x-m$, and a
second time in the term $r(x) = rx/(x+n)$ which can be interpreted as an approximation of a population dynamics where
the differences between fertile and non-fertile are not explicitly modelled.
The one with a unique Allee effect is given by
\begin{equation}\label{Alleesimple}
	\begin{split}
		\dot{x}&=\left(r\left( 1-\frac{x}{K}\right)(x-m)-qy\right)x, \\
		\dot{y}&= (px-c)y.
	\end{split}
\end{equation}
We recall that this last system is a particular case of \eqref{Alleefg}.

For system \eqref{Alleedoble} there exist always three equilibria, two of them in the boundary. In the case with strong Allee effect there exists another boundary equilibrium.

The authors prove the existence of a subset of parameter values for which there are two limit cycles when the Allee effect is either strong or weak in system \eqref{Alleedoble}. This is a different result by comparing with system \eqref{Alleesimple}, generated by the double Allee effect, because  system \eqref{Alleesimple} has a unique limit cycle generated by Hopf bifurcation, surrounding the unique positive equilibrium.

In system \eqref{Alleedoble} there exists a separatrix curve determined by the unstable manifold of equilibrium point $(m, 0)$ or the origin. Then, there are trajectories close this separatrix with different $\omega$-limit sets for the same values of the parameters, so the system is highly sensitive to initial conditions. For a same set of parameters different behaviours are possible: the extinction of two populations, the coexistence, and the oscillation of both populations. This is in acordance with the results in \cite{Olivares2010}.

\begin{figure}[h]
	\centering
	\subfloat[Coexistence in system \eqref{Alleedoble}.]{\includegraphics[width=55mm]{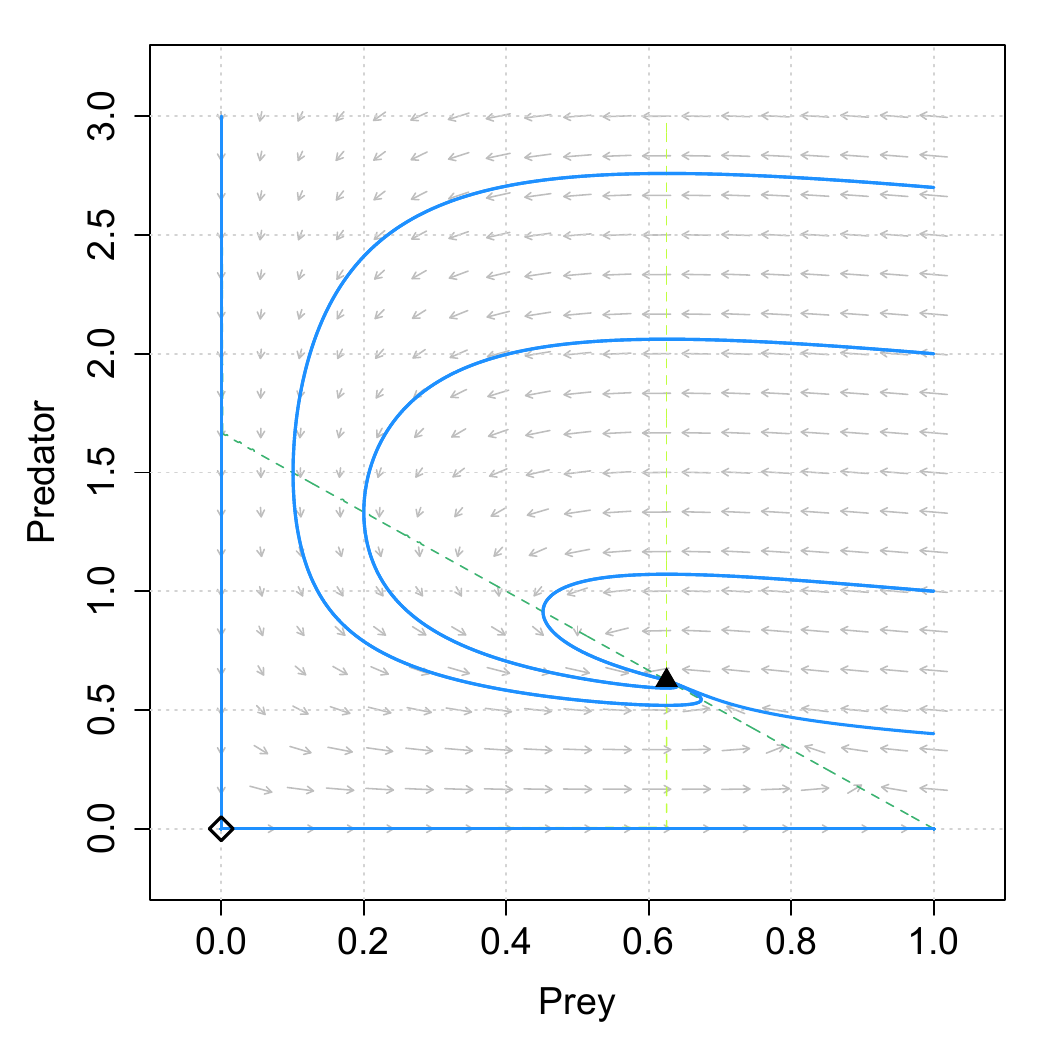}}
	\subfloat[Coexistence in system \eqref{Alleesimple}.]{\includegraphics[width=55mm]{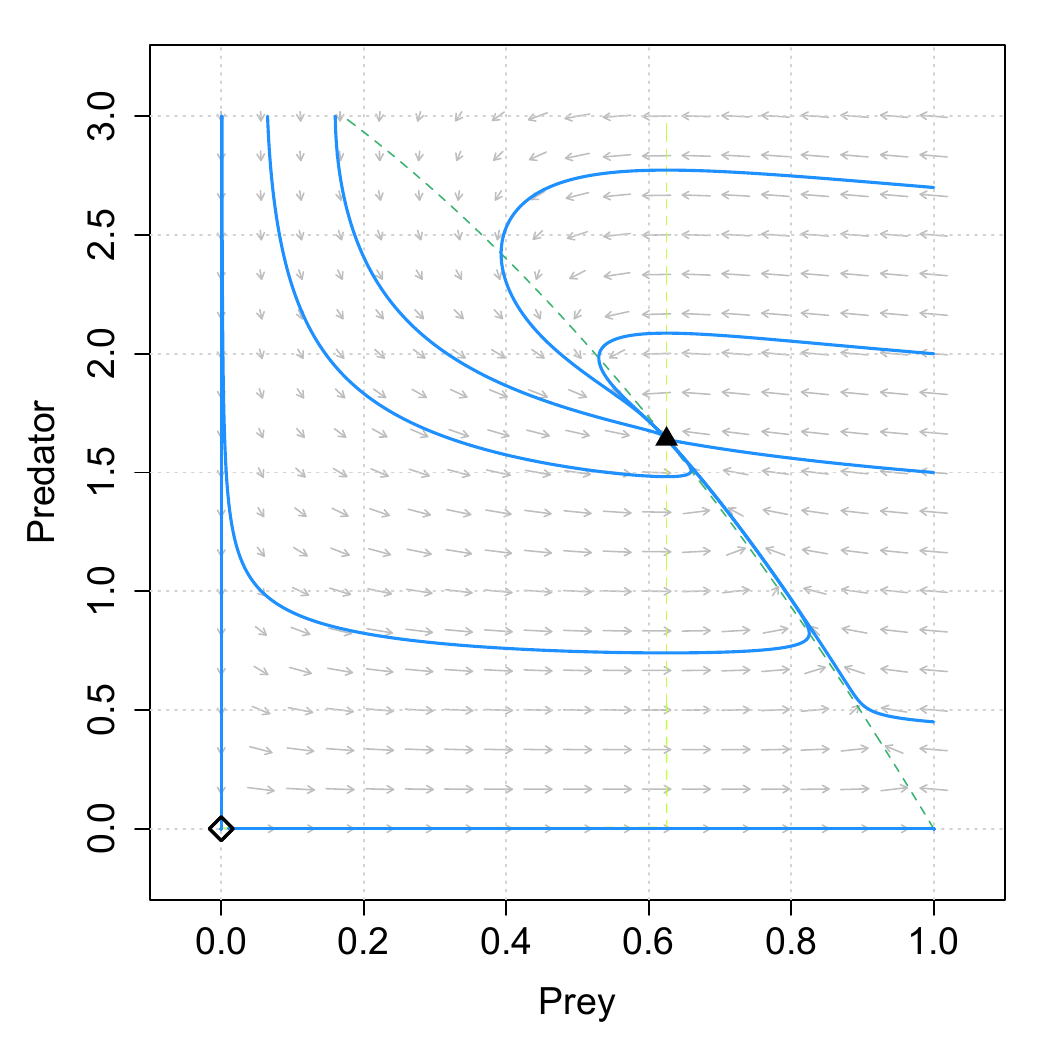}}\\
	\caption{Dynamic behaviours with $r=2$, $n=2$, $K=1$, $m=-2$, $q=1.2$, $p=0.8$ and $c=0.5$. }
	\label{fig:(13)(14)a}
\end{figure}

The separatrix curve disappears when $m<0$ and in this case the origin is a saddle point in both systems \eqref{Alleedoble} and \eqref{Alleesimple}. We give an example in Figure \ref{fig:(13)(14)a} of how the systems \eqref{Alleedoble} and \eqref{Alleesimple} show an analogous behaviour but in the system with double Allee effect the predator density in the coexistence equilibrium is smaller than in the case with simple Allee effect.

\begin{figure}[H]
	\centering
	\subfloat[Phase portrait for system \eqref{Alleedoble}.]{\includegraphics[width=55mm]{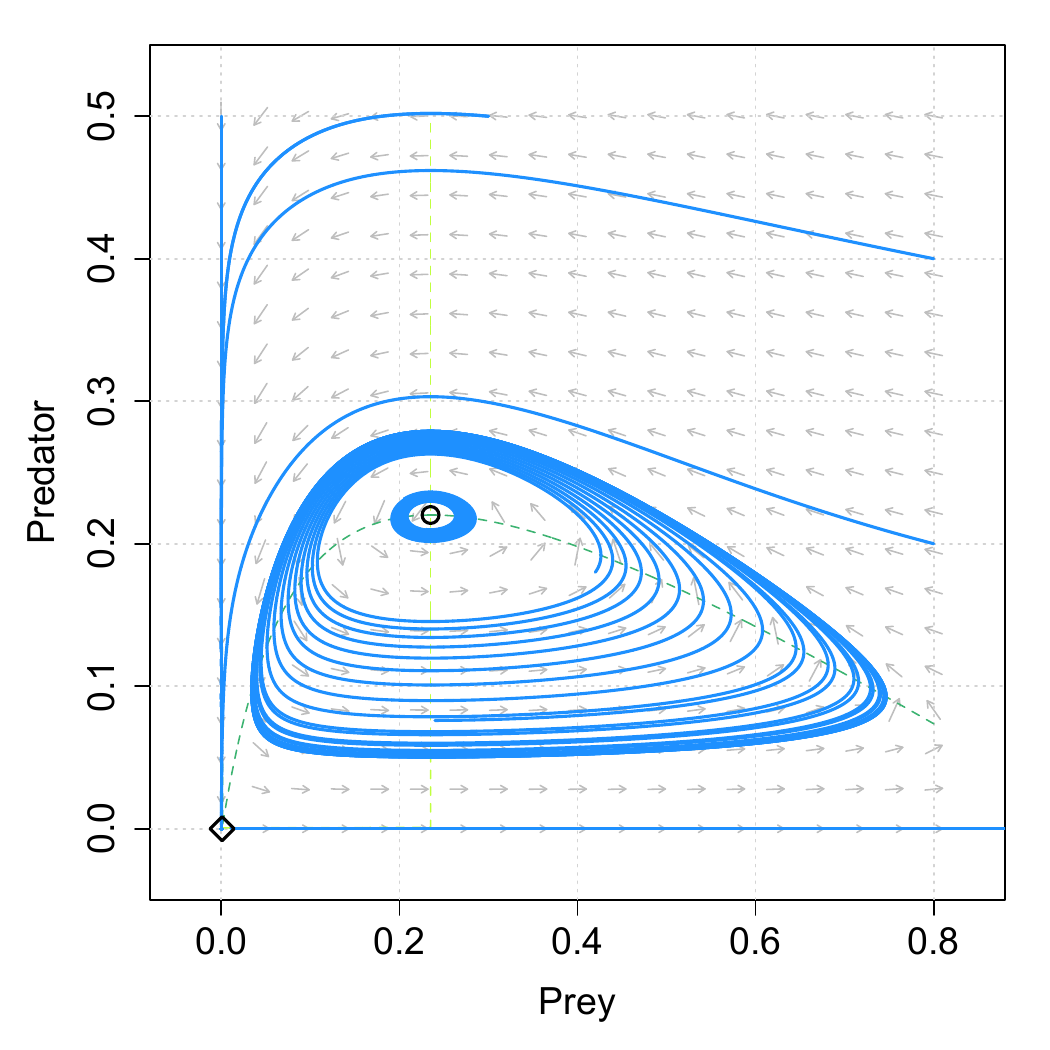}}
	\subfloat[Phase portrait for system \eqref{Alleesimple}.]{\includegraphics[width=55mm]{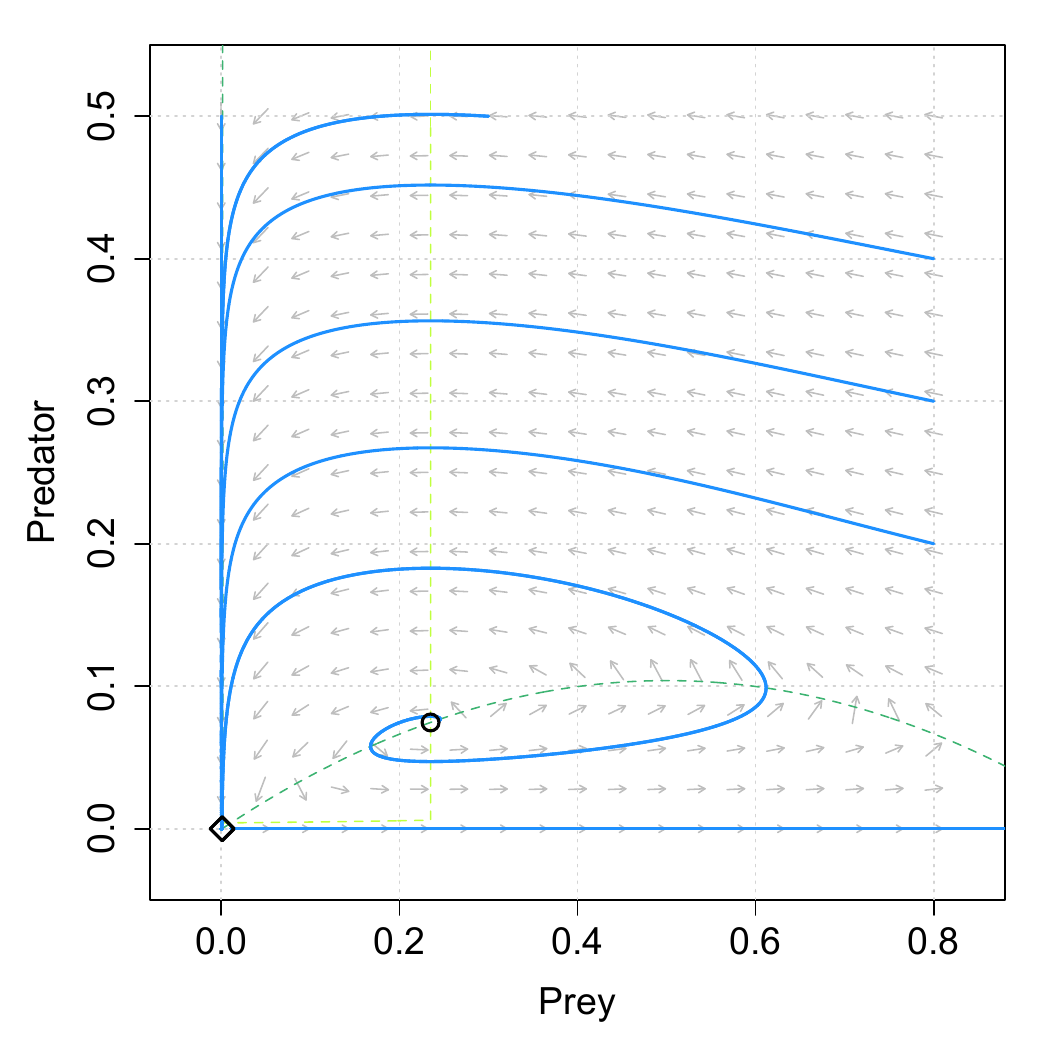}}\\
	\caption{Dynamic behaviours with $r=0.5$, $n=0.10375$, $K=1$, $m=0.001$, $q=1.2$, $p=0.175$ and $c=0.041125$. }
	\label{fig:(13)(14)b}
\end{figure}

In Figure 	\ref{fig:(13)(14)b} we compare the phase portraits of systems \eqref{Alleedoble} and \eqref{Alleesimple} with the same values of the common parameters. We note that for the system \eqref{Alleedoble} there exist two limit cycle surrounding the positive equilibrium, while it does not occur for the system with simple Allee effect.

\vspace{0.5cm}
With similar techniques, the influence of the Allee effect is studied in many other works, combined with different particular biological conditions. 
For example, the influence of the inclusion of different characteristics in Gause type predation models has been studied. In particular, in \cite{Olivares2019} the competition among predators is included, and in \cite{Olivares2011NA,Olivares2011NAcorr} the hyperbolic functional response, which is a particular case of Holling type II functional response, is considered. The authors study both the strong and weak Allee effect. A generalization of this last work is given in \cite{Olivares2013}, with different functional responses. In this last one, the results show a marked difference between cases of strong and weak Allee effect, in constrast with \cite{Olivares2019} where there exist minor differences among both cases.

\vspace{0.5cm}

Instead of assuming a functional response of type I, as in the previous models, now we are presenting a model with type II functional response. This is the case of the work of D. Sen, S. Petrovskii, S. Ghorai and M. Banerjee, \cite{Allee_genomn}. They also consider generalist predators, i. e.  the ones that feed not only on the species of prey. The predator species is subject to an Allee effect, and this model is motivated by the real world example of the phytoplankton-zooplankton system. While the zooplankton growth is known to be subject to the Allee effect, there is considerable evidence that the phytoplankton growth is logistic \cite{plankton}. 

They show that increasing the carrying capacity of preys helps the population to survive in a coexistence equilibrium. Their model is capable of producing bi-stable dynamics for a reasonable range of parameter values. The obtained dynamics is quite rich and there exist different types of bifurcation. The system they consider is:
\begin{equation}\label{15}
	\begin{split}
		\dot{x}&= r_1 x \left( 1-\frac{x}{k_1}\right) -\frac{axy}{1+ahx},\\[0.1cm]
		\dot{y}&=r_2y\left(1-\frac{y}{k_2} \right) (y-m) + \frac{e a x y }{1+ah x},
	\end{split}
\end{equation}
which has a Holling Type II functional response. Taking $x=N/k_1$ and $y=P\sqrt{(r_2/(k_2r_1))}$ they get a dimensionless system:
\begin{equation}\label{genP}
	\begin{split}
		\dot{x}&=x(1-x)-\frac{\alpha x y}{\delta+x},\\
		\dot{y}&=y(\eta-y)(y-l)+\frac{\alpha_1 x y}{\delta+x},
	\end{split}
\end{equation}
where $\alpha=(1/(h k_1))\sqrt{k_2/(r_1r_2)}$, $\delta=1/(ah k_1)$, $\eta=\sqrt{r_2 k_2/r_1}$, $l=m\sqrt{r_2/(k_2r_1)}$ and $\alpha_1=e/(r_1h)$.

This system has the trivial equilibrium $E_0=(0,0)$ and three boundary equilibria $E_1=(1,0)$, $E_2=(0,l)$ and $E_3=(0,\eta)$. The system can also present up to four positive equilibria $E_j^*$ obtained as the intersection points of non-trivial prey and predator nullclines in the first quadrant:
\begin{equation}
	\begin{split}
		&y-\frac{1}{\alpha}(1-x)(\delta+x)=0,\\
		&(\eta-y)(y-l)+\frac{\alpha_1x}{\delta+x}=0.
	\end{split}
\end{equation}
Obviously, the number of interior equilibria depends on the value of the parameters, shape and position of the nullclines and position of the points of intersection of the nullclines on the positive axes with respect to the boundary equilibria.

Stability conditions for all boundary equilibria are obtained by applying standard  linear stability analysis, but for the positive equilibria the authors make a graphical approach, since the expressions of $E_j^*$ are not available explicitly.

For the boundary equilibria they prove that $E_0$ is always a saddle; $E_1$ is a saddle if $\alpha_1>(1+\delta)\eta l$ and locally asymptotically stable if $\alpha_1<(1+\delta)\eta l$; 
$E_2$ is a saddle if $\delta/\alpha<l$ and an unstable node if $\delta/\alpha>l$; $E_3$ is locally asymptotically stable if $\delta/\alpha<\eta$ and it is a saddle if $\delta/\alpha>\eta$. They also give a description of the basin of attractions of the equilibria.

Regarding bifurcations appearing on the system, the authors proof the existence of saddle-node and transcritical bifurcations through which the positive equilibria appear or dissapear and the stability of the boundary equilibria changes. Their results are the following: The point $E_1$ undergoes a transcritical bifurcation when $\delta$ passes the threshold $\delta_{TC1}=\alpha_1/(\eta l)-1$ for $\alpha_1>\eta l$. $E_1$ undergoes another transcritical bifurcation at the threshols $\delta_{TC2}=\alpha l$. At last, $E_3$ undergoes a transcritical bifurcation at $\delta_{TC3}=\alpha\eta$.

They also prove analytically the existence of saddle-node bifurcations, more precisely: the system \eqref{genP} undergoes a saddle-node bifurcation when the curve $y=p(x)$ touches 
$y=q_1(x)$ at $\overline{E}(\overline{x},\overline{y})$ such that $\overline{x}$ is a double root of $h(x)=0$ with $\overline{x}\in(0,1)$ and $\overline{y}\in(0,l)$. System \eqref{genP} undergoes another saddle-node bifurcation at $\alpha=\alpha_{SN2}$ when $y=p(x)$ touches $y=q_2(x)$ at $\overline{\overline{E}}(\overline{\overline{x}},\overline{\overline{y}})$ such that $\overline{\overline{x}}$ is a double root of $h(x)=0$ and $\overline{\overline{y}}=p(\overline{\overline{x}})=q_2(\overline{\overline{x}})$ with $0<\overline{x}<1$. The expressions for $p(x), q_1(x), q_2(x)$ and $h(x)$ are:
\begin{equation}
	\begin{split}
		p(x)&=\frac{1}{\alpha}(1-x)(\delta+x),\\
		q_1(x)&=\frac{1}{2}\left( (\eta+l)-\sqrt{(\eta-l)^2+\frac{4\alpha_1x}{x+\delta}}\right), \\
		q_2(x)&=\frac{1}{2}\left( (\eta+l)+\sqrt{(\eta-l)^2+\frac{4\alpha_1x}{x+\delta}}\right),\\
		h(x)&=x^5+(3\delta-2)x^4+ (3\delta^2+6\delta+1+\alpha(\eta+l))x^3+(\delta^3-6\delta^2+\delta(3+2\alpha(\eta+l))\\
		&-\alpha(\eta+l))x^2+(-2\delta^3+\delta^2(3+\alpha(\eta+l))-2\alpha\delta(\eta+l)-\alpha^2(\alpha_1-\eta l))x\\
		&+(\delta^3-\alpha \delta^2(\eta+l)+\alpha^2\delta\eta l).
	\end{split}
\end{equation}

As the expressions of the positive equilibria are difficult to find, the authors can not give the analytical conditions for Hopf-bifurcation, but they find the numerical threshold value of the parameter $\alpha$ at which the system \eqref{genP} undergoes this kind of bifurcation. For this they use the software Maple. They prove numerically that the first Lyapunov coefficient is negative, so the bifurcating limit cycle is stable and the bifurcation is supercritical.

They also study numerically the so called Bogdanov-Takens bifurcation. It is a  local co-dimension 2 bifurcation which occurs at the point where the bifurcation curves corresponding to the saddle-node bifurcation and the Hopf bifurcation intersect when we 
draw them togehter in a two parametric bifurcation plane. The authors give the threshold value of this bifurcation.

The study of the complete dynamics of the system \eqref{genP}  carried out in \cite{Allee_genomn}  shows with a well-founded analysis that the dynamics of the predator-prey system can be very complex also in systems with only two species. The authors show that the dynamics of the generalist predator-prey system assuming that the predator is subjected to Allee effect is much more complicated than in the case of the specialist predator. The prey-predator system with a generalist predator has more sensible biological properties,  because the extinction of prey does not lead to the extinction of predator.
While the specialist predator system has some unrealistic properties, the generalist predator system apparently does not. One remarkable property of the specialist case is its counter-intuitive response to the system enrichment described by an increase in either its linear per capita growth rate or its carrying capacity. 
The complexity is also increased because the inclusion of the Allee effect in the predator causes different types of bifucation.

The work done in \cite{Allee_genomn} leaves some open questions as if the complexity of the system may increase considering Holling type III or  ratio-dependent functional response instead of the Holling type II functional response. 

\color{black}

\section{The influence of fear effect}\label{sec:fear}

In a biological system where predators and prey coexist, the behaviour of each of these species can affect both populations. 
The way of life of the prey affects the evolution of the ecosystem, for example when prey, aware of the presence of their predators, feel fear and act accordingly making hunting more difficult for the predators. This fact drives prey to avoid direct predation, which can increase short-term survival of prey, but may cause a long-term decrease in prey population as a consequence \cite{Cresswell}.

There exist experiments that support the theoretical reasonings about the effect that anti-predator behaviors may have, as the ones mentioned in the introduction \cite{songbirds,creel,wirsing}. All these experiments on birds and vertebrate species reported the same kind of conclusions: even though there is no direct killing between predators and prey, the presence of predators cause a reduction in prey population due to anti-predator behaviors.

We analyze different models with the presence of fear effect. First of all, in \cite{fearWang} the influence of fear effect in two different systems is studied. This two systems consider respectively a Holling type I functional response and a Holling type II functional response. In this model predators feed only on prey, so in \cite{fear} this model is generalized into one in which predators are omnivorous. Then we present a study of the impact of fear on a Leslie-Gower model with hunting cooperation \cite{miedocaza}. As we said in the introduction, the combination of the fear effect and hunting cooperation is supported by many ecological evidences.
Next we present two examples where the fear effect is combined with Allee effect, the first one in \cite{Sasmal} which has a multiplicative Allee effect and the second in \cite{fear-alleeaditive} with an additive Allee effect.
Finally, in \cite{fear-preyrefuge} the influence of prey refuge is studied, including it on a model with Holling type II functional response. The prey refuge is also considered in \cite{pr-fear-allee}, where it is included in the model proposed in \cite{Sasmal} and combined with an Allee effect.

\vspace{0.5cm}
We start with the work of X. Wang, L. Zannete and X. Zou in \cite{fearWang} where they propose the following model,
\begin{equation}\label{siswang}
	\begin{split}
		\dot{x}&=r x f(s,y)- dx - ax^2- g(x)y,\\
		\dot{y}&= y(-n + c g(x)), 
	\end{split}
\end{equation}
which modifies a classic predator-prey model multiplying the production term by a factor $f(s,y)$ related with the cost of anti-predator defence due to fear. The parameter $s$ represents the level of fear and by the biological meaning the following conditions are asumed: 
\begin{equation*}
	\begin{split}
		&f(0,y)=1, \;\;\;\;  f(s,0)=1, \;\;\;\;  \lim_{s\to \infty}f(s,y)=0, \;\;\;\; \lim_{y\to \infty}f(s,y)=0,\\
		& \frac{\partial f(s,y)}{\partial s}<0, \;\;\;\; \frac{\partial f(s,y)}{\partial y}<0.
	\end{split}
\end{equation*}

In \cite{fearWang} two types of functional responses are considered: a Holling type I functional response or lineal response, $g(x)=x p$, and a Holling type II functional response, $g(x)~=xp/(1+qx)$. 

It is assumed that predators feed exclusively on this species of prey. Consequently, only two boundary equilibria appear, since there could never be a population of predators in the absence of prey, and therefore the equilibrium on the vertical axis which appears on the classical model does not exist in this case. 

The authors prove that if $r<d$ then both populations get extinct, independently of the fear effect and the particular predation mechanism, because $x'(t)\leq (r-d)x$, and by a comparison argument and the theory of asymptotically autonomous systems \cite{Castillo-Chavez}, they obtain that $x(t)$ and $y(t)$ tend to zero when $t$ tends to infinity. Therefore, it is enough to consider the case when $r > d$.

Taking the linear functional response,  if $r\in (d,d+an/cp)$ then there is not 
positive equilibrium, i. e, there is not coexistence of both species, and the equilibrium on the boundary, $E_1=((r-d)/a,0)$ is globally asymptotically stable, so predators get extinct. If $r>d+an/cp$ there exists also a positive equilibrium, $E^*$, which is globaly asymptotically stable.  We note that the threshold value $d+an/cp$ does not depend on the parameter $s$, so it seems that fear effect has no influence on the dynamics in this case. 

With the Holling type II functional response, they consider that the term related with the fear effect is
\begin{equation*}
	f(s,y)=\frac{1}{1+sy}.
\end{equation*}
If $cp\leq nq$, then there is not positive equilibrium and the boundary equilibrium $E_1$ is globally asymptotically stable. If $cp>nq$, applying a transformation on the variables $t,x$ and $y$ they get the Kolmogorov system

\begin{equation}
	\begin{split}
		\dot{x}&= x (a_1+a_2 x - a_3 y - a_4 x y - a_5 x^2 - a_6 y^2 - a_5 x^2y),\\
		\dot{y}&=y(x-1)(1+y), 
	\end{split}
\end{equation}
where $a_i>0$ for $i=1,\dots,6$, whose analysis allows them to determine the stability of the equilibrium points of the original system, obtaining the following results. 
If $(r-d)(cp-nq)>an$ there exist only two equilibrium points, the origin which is unstable and $E_1$ which is stable.
If $(r-d)(cp-nq)>an$ the origin and $E_1$ are unstable and there exist a positive equilibrium $E^*$ which stability depends on the parameters. $E^*$ is locally asymptotically stable if 
\begin{equation}\label{4.19}
	r >\frac{an}{cp-nq}+d \;\; \text{ and } \;\; 
	r \leq d + \frac{a(cp + nq )}{q (cp - nq)}, 
\end{equation}
or 
\begin{equation}
	r> d + \frac{a(cp + nq)}{q(cp - nq)} \;\; \text{ and } \;\; 
	s >\frac{q (cp - nq)^2((r - d )q(cp - nq)-a(cp+nq))}{c^2 p a (q d (cp-nq)+a(cp + nq))}, 
\end{equation}
and $E^*$ is unstable if
\begin{equation}
	r> d + \frac{a(cp + nq)}{q(cp - nq)} \;\; \text{ and } \;\; 
	s <\frac{q (cp - nq)^2((r_0 - d )q(cp - nq)-a(cp+nq))}{c^2 p a (q d (cp-nq)+a(cp + nq))}.
\end{equation}

When the birth rate of prey is not large enough to support oscillations, there exists a globally asymptotically stable positive equilibrium, and when the birth rate is large enough to support oscillations, the positive equilibrium can be also locally asymptotically stable if prey are sensitive enough to perceive potential dangers and show an anti-predator behaviour motivated by fear. So the conclusion is that in some cases, the fear effect can stabilize the system.

It is proved that the asymptotical stability of $E^*$ is global if \eqref{4.19} holds and also 
\begin{equation}
	1 \leq \frac{r q}{cp - nq}.
\end{equation}

In addition, in the case with Holling type II functional response, it is possible the existence of limit cycles. In this case, the positive equilibrium becomes unstable and a limit cycle appears as a result of a Hopf bifurcation. The impact of the fear effect on this Hopf bifurcation is also studied in \cite{fearWang}.

Although the authors choose $f(s,y)=1/(1+sy)$ for the theoretical analysis of the model, they also conduct some numerical simulations with another decreasing functions, in particular: 
\begin{equation*}
	f(s,y)=e^{-sy} \;\; \text{ and } \;\; 
	f(s,y)=\frac{1}{1+s_1y+s_2y^2}.
\end{equation*} 

The behaviour observed in the simulations regarding the existence of Hopf bifurcation with these other functions is similar to that obtained for the first chosen function. In general, the same results are obtained for a general monotone decreasing function of $y$.

The model analysed in this paper assumes that perceived predation risks only reduce the birth rate and survival of offspring, and does not include the possible impact on the mortality rate of adult prey. They also do not include the case where fear affects intraspecific competition. The reason for not incorporating them in the models is that there is not yet experimental evidence for both effects. However, in \cite{songbirds} and \cite{clinchy} they argue that fear can increase the mortality rate of adults due to long-term physiological impacts. There are also theoretical arguments in \cite{Cresswell} that the effect of fear may change the strength of intraspecific competition due to the complexity of the food web.
\color{black}

\vspace{0.5cm}
In the previous model, predators feed only on the prey species, but recently,
Z. Zhu, R. Wu, L. Lai and X. Yu in \cite{fear} have proposed a model where predators are omnivorous. As we saw in Section \ref{sec:Allee} for the case of populations with Allee effect, the consideration of omnivorous or generalist predators affects the dynamics of the system.
Here the authors assume predators can feed on other species, particularly when the main prey species becomes extinct. This characteristic is incorporated on the system with lineal functional response considered in \cite{fearWang}, by imposing a logistic growth of predators in absence of prey  obtaining:
\begin{equation}\label{25}
	\begin{split}
		\dot{x}&=\frac{r x }{1+sy}- dx - ax^2 - p x y, \\
		\dot{y}&=c p x y + ny - d_1 y^2.
	\end{split}
\end{equation}

Depending on the values of the parameters, the system can exhibit four equilibria. Two of them are the boundary equilibriums that also appeared on the system without omnivorous predators, $E_0$ and $E_1$, and now it appears also a new boundary equilibrium $E_2=(0,n/d_1)$, because  predators can survive in absence of prey. Here, the equilibrium points $E_0$ and $E_1$ are unstable, so the dynamic behaviour of both systems are quite different. Furthermore, the equilibrium $E_2$ is globally asymptotically stable under certain conditions  which shows that the predator species can be permanent despite the extinction of the prey species. When the positive equilibrium exists, it is globally asymptotically stable, as in the case with no omnivorous predators.

\begin{figure}[h]
	\centering
	\subfloat[Solutions for $x(0)=0.1$ and $y(0)=0.1$]{\includegraphics[width=65mm]{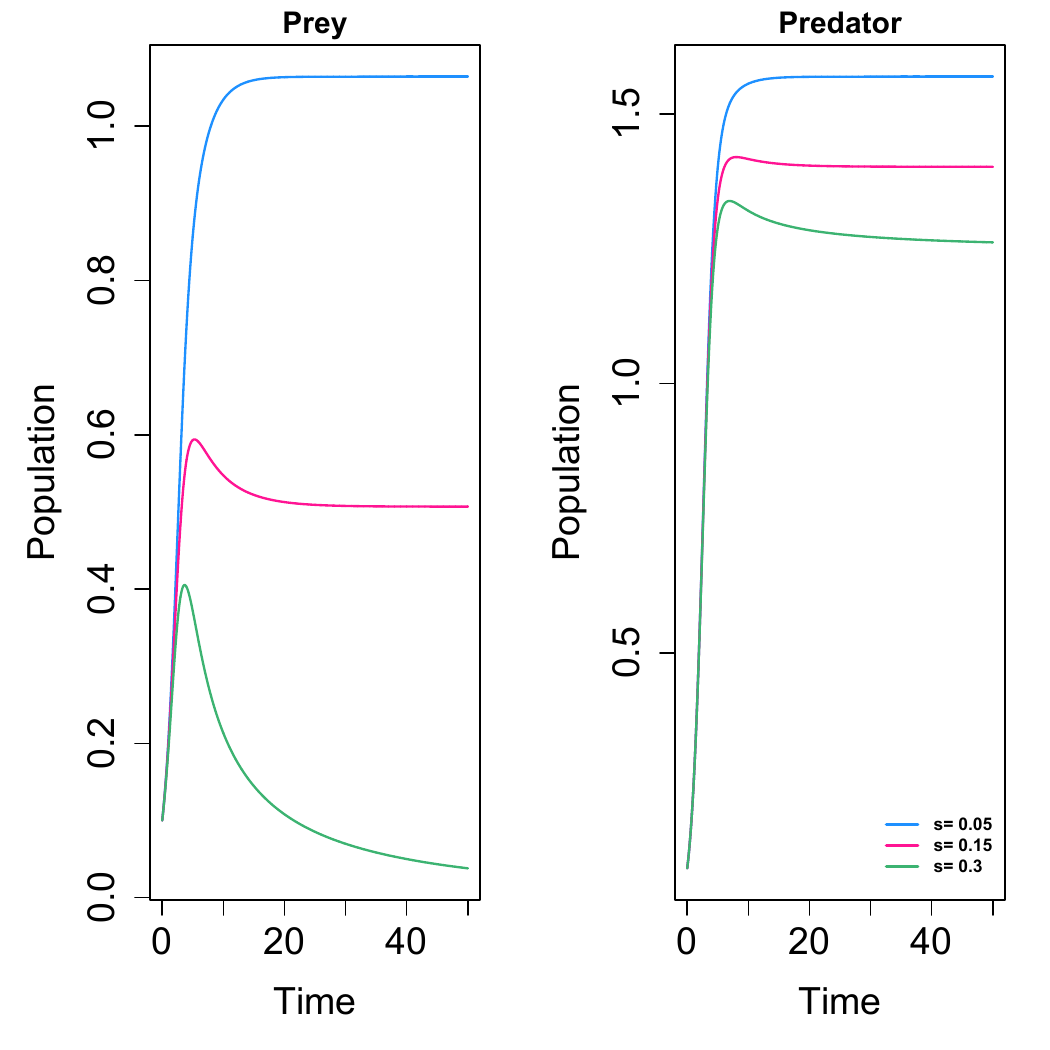}}\\
	\subfloat[Variation of $x^*$]{\raisebox{0.4cm}{\includegraphics[width=40mm]{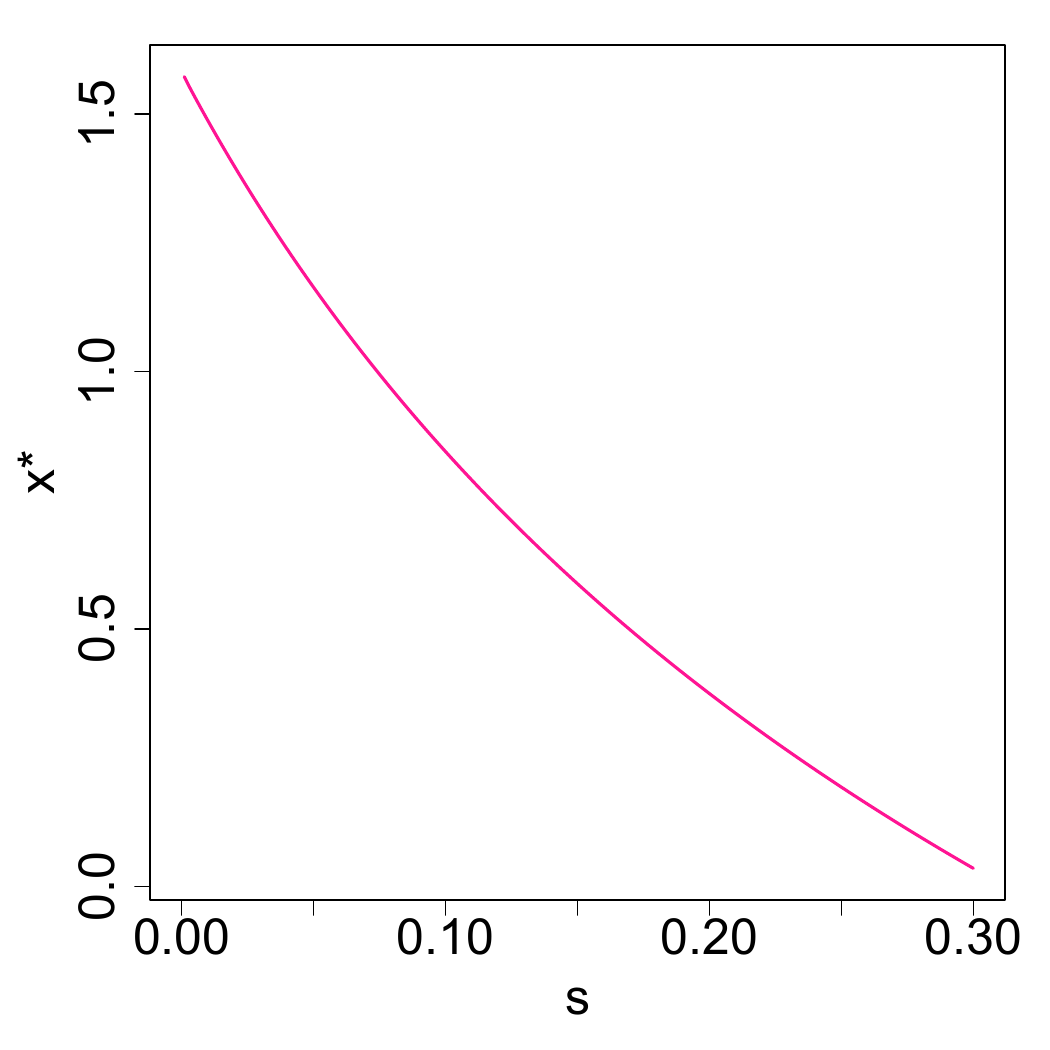}}}
	\subfloat[Variation of $y^*$]{\raisebox{0.4cm}{\includegraphics[width=40mm]{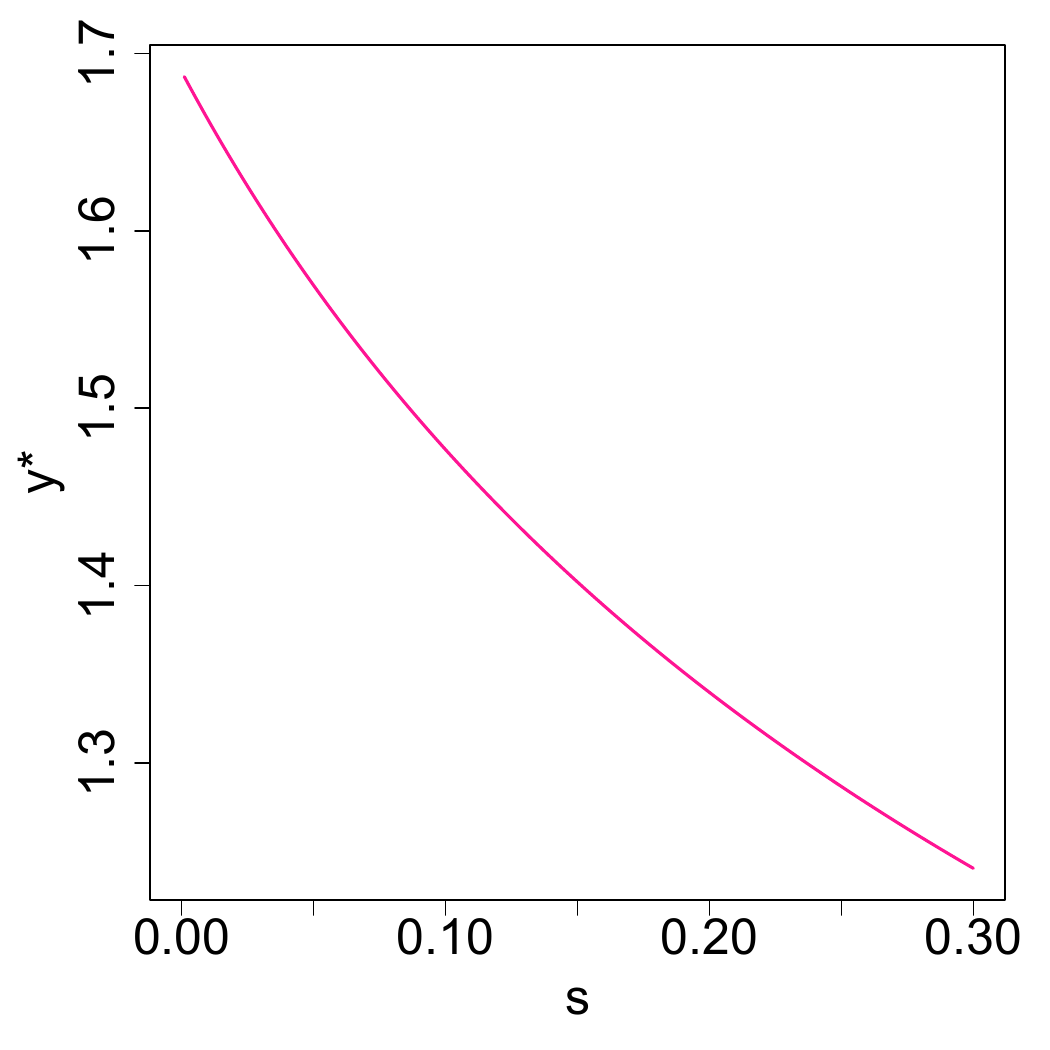}}} 
	\caption{Variation of the positive equilibrium and the solutions in function of parameter $s$. The considered values of the parameters are $r = 1.6$, $d=0.8$, $a=0.2$, $p=0.3$, $c=0.8$, $n=1$ and $d_1=0.8$. }
	\label{fig:(25)k}
\end{figure}

In addition, the authors conclude that the density of both prey and predator species decrease when the fear parameter $s$ increases. In Figure \ref{fig:(25)k} 
we show an example of how the variation of the fear parameter affects the solutions and the final population density of each species in the equilibrium point. 

\begin{figure}[H]
	\centering
	\subfloat[Extinction of both species in system \eqref{siswang}]{\includegraphics[width=55mm]{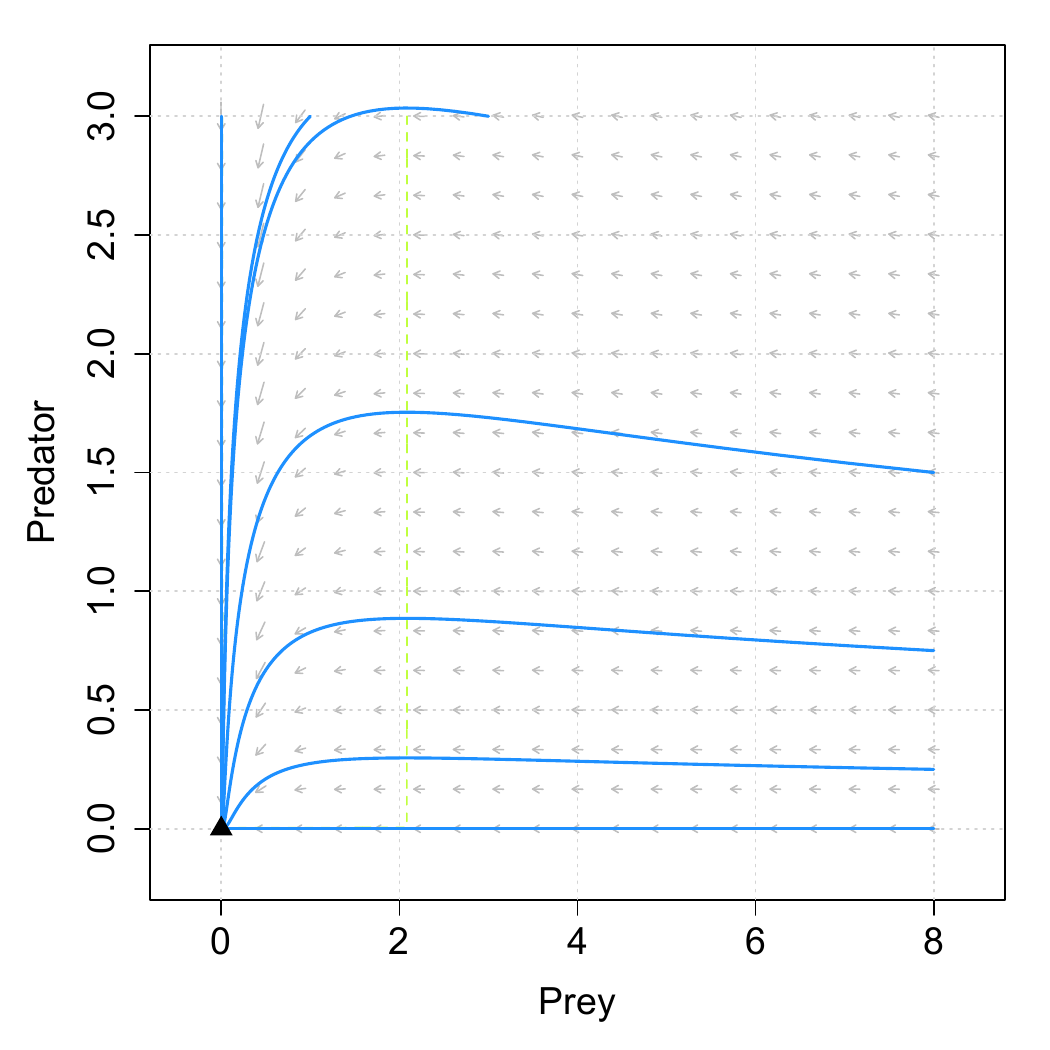}}
	\subfloat[Extinction of prey in system \eqref{25}]{\includegraphics[width=55mm]{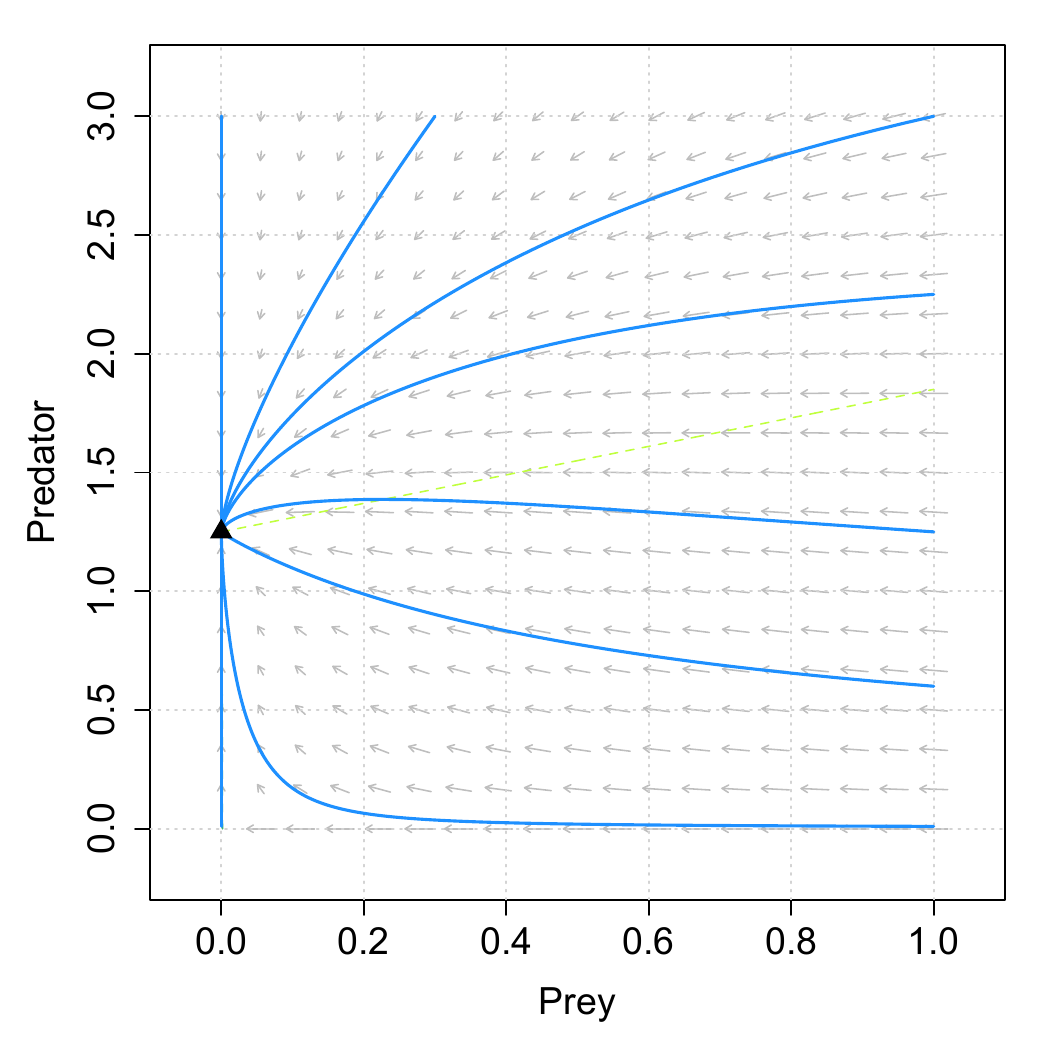}}\\
	\caption{Dynamic behaviours with $r=0.5$, $s=3$, $d=0.8$, $a=1.5$, $p=0.6$, $c=0.8$, $d_1=0.8$ and $n=1$.}
	\label{fig:(25)(19)a}
\end{figure}

\begin{figure}[H]
	\centering
	\subfloat[Extinction of predators in system \eqref{siswang}]{\includegraphics[width=55mm]{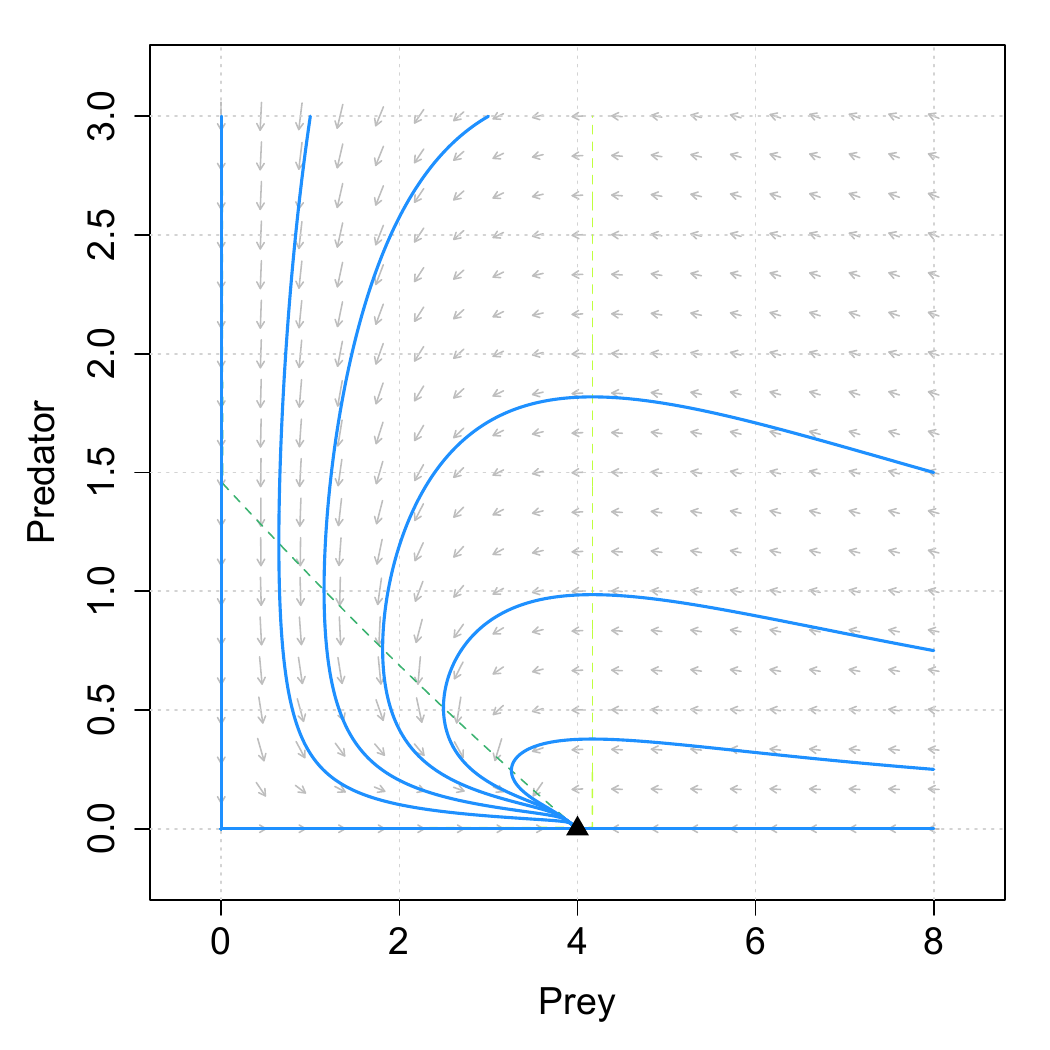}}
	\subfloat[Coexistence in system \eqref{25}]{\includegraphics[width=55mm]{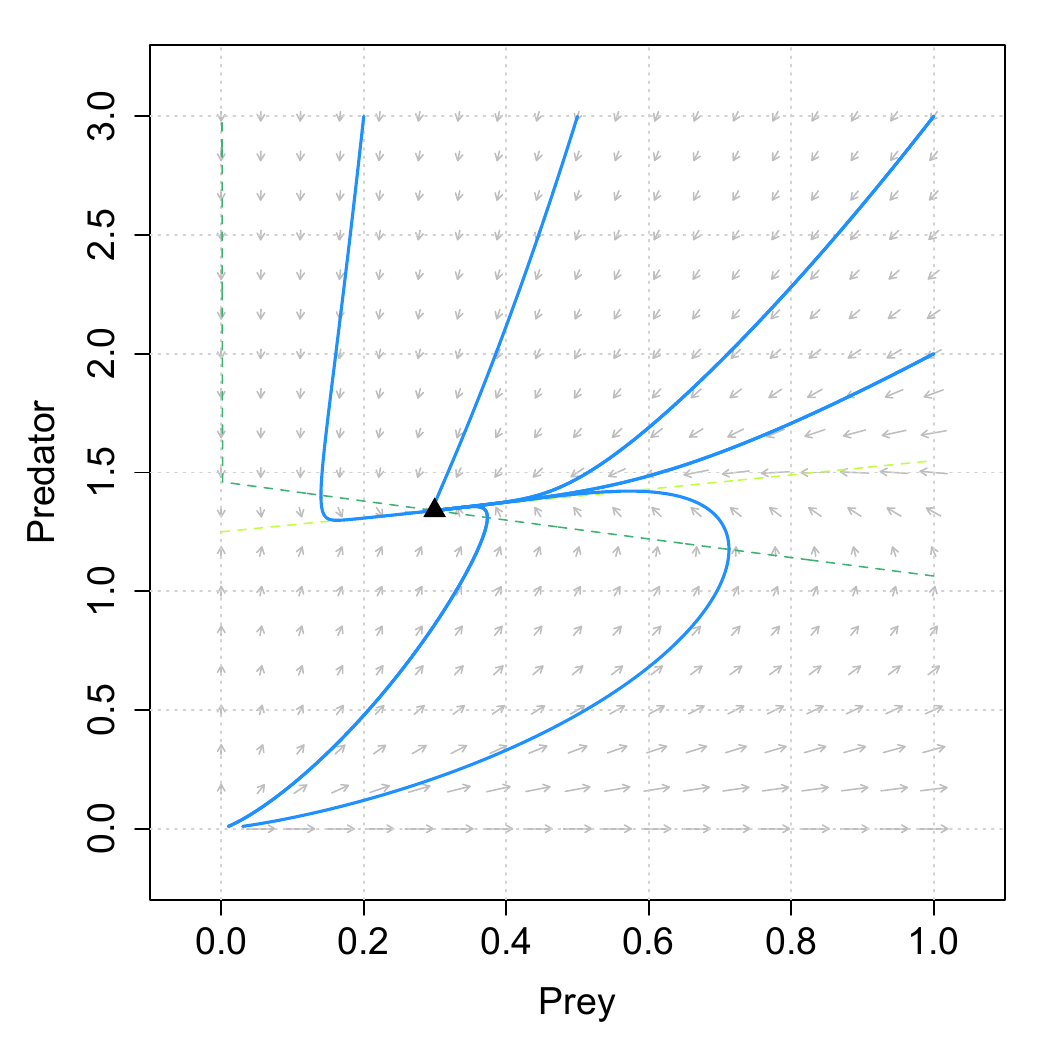}}
	\caption{Dynamic behaviours with $r=1.6$, $s=0.2$, $d=0.8$, $a=0.2$, $p=0.3$, $c=0.8$, $d_1=0.8$ and $n=1$.}
	\label{fig:(25)(19)b}
\end{figure}

It is concluded that fear effect can cause the extinction of the species, in opposition to the conclusions of \cite{fearWang}, where the fear did not affect the dynamics. In Figures \ref{fig:(25)(19)a} and \ref{fig:(25)(19)b} we give examples of the different behaviours obtained for systems \eqref{25} and \eqref{siswang} (with linear response) with the same values of the common parameters, so it is clear that the inclusion of the omnivorism affects the dynamics.

\vspace{0.5cm}
Considering the importance of the predator predation strategy and the prey avoidance strategy it is convenient to analyze how both prey and predators adjust their behaviours to obtain maximum benefits and increase their biomass for each of them. In \cite{miedocaza},  S. Pal, N. Pal, S. Samantha and J. Chattopadhyay consider a modified Leslie-Gower predator-prey model where predators cooperate during hunting and due to fear of predation risk prey populations show anti-predator behaviour.  The impact of hunting cooperation and fear effect on the dynamics of the system are analyzed.

The authors consider the following system with Holling functional response of type II, fear effect and Berec's encounter-driven functional response
\cite{Berec}:

\begin{equation}\label{fearcaza}
	\begin{split}
		\dot{x}&=x\left( \frac{b}{1+\beta y}-d-x-\frac{(1+\alpha y)y}{p+(1+\alpha y)x}\right),\\[0,1cm]
		\dot{y}&=y \left( 1- \frac{q y}{r + (1+\alpha y)x}\right).
	\end{split}
\end{equation}
Note that in absence of hunting cooperation $(p=0)$, the previous model coincides with the modified Leslie-Gower model.

It is proved that solutions of system \eqref{fearcaza} starting at positive conditions are bounded if $b>d$ and $q>(b-d+\varepsilon_1)$ with $\varepsilon_1>0$. Moreover, system \eqref{fearcaza} is permanent if $b/(1+\beta M_2)-d - M_2(1+\alpha M_2)/p >0$, where $M_2$ is the upper bound of predator population. 

The model has the following boundary equilibria: $E_0=(0,0)$, $E_1=(b-d,0)$ and $E_2=(0,r/q)$. It can present also a positive equilibrium $E^*=(x^*,y^*)$. The equilibium points $E_0$ and $E_1$ are both unstable. The equilibrium $E_2$ is locally asymptotically stable if the cooperation stength $\alpha$ is greater than the threshold value
\begin{equation*}
	\frac{b p q^3}{r^2(\beta r + q)}- \frac{d p q^2}{r^2}- \frac{q}{r},
\end{equation*}
or if the fear strength $\beta$ is greater than the threshold value 
\begin{equation*}
	\frac{b p q^3}{r^2(\alpha r + q)+d p q^2 r}- \frac{q}{r}.
\end{equation*}

By using Routh-Hurwitz criterion, the positive equilibriun $E^*$ is locally asymptotically stable if 
\begin{equation*}
	\begin{split}
		T(x^*,y^*)&=-1-x^*+\frac{\alpha x^*}{q}+\frac{x^* y^* (1+\alpha y^*)^2}{(q y^* + p -r)^2}<0  \;\; \text{ and}\\[0.1cm]
		D(x^*,y^*)&=\left( \frac{\alpha x^*}{q}-1\right) \left( -x^*+\frac{x^* y^* (1-\alpha y^*)^2}{(q y^* + p -r)^2}\right)\\[0.1cm]
		&+ \frac{1+\alpha y^*}{q}\left( \frac{b\beta x^*}{(1+\beta y^*)^2} + \frac{(qy^*-r)(p+x^*+\alpha x^*y^*)+p \alpha x^* y^*}{(q y^* + p - r)^2}\right) >0.
	\end{split}
\end{equation*}

Regarding the existence of bifurcations, this system presents a Hopf bifurcation at $\alpha=\alpha_h$ from the equilibrium point $E^*$ if and only if $D(x^*(\alpha),y^*(\alpha))>0$ and $\frac{d}{d\alpha}T(x^*(\alpha),y^*(\alpha))\neq0$ in $\alpha=\alpha_h$.

In order to obtain the stability and direction of the Hopf equilibrium, the first Lyapunov coefficient $l$ is calculated, applying the results on \cite{Wiggins}. They obtain that this bifurcation can be supercritical or subcritical depending on the parameters. 

We note that when $l=0$ the system exhibits a generalized Hopf-bifurcation, which is also called Bautin bifurcation, at which the interior equilibrium has a pair of purely imaginary eigenvalues. The Bautin bifurcation point separates branches of subcritical an supercritical Hopf-bifurcation in the parameter plane. 

The Bogdanov-Takens bifurcation from the equilibrium point $E^*$ is also analyzed. The employed techniques are given in \cite{kuznetsov}.

Some of the important conclusions of this work are the following. 
In the absence of fear, if we increase the strength of hunting cooperation, then the system becomes unstable around the positive equilibrium and undergoes multiple Hopf-bifurcations, the first of them is supercritical and the second one is subcritical. The coexisting equilibrium can lose the stability via supercritical Hopf bifurcation with an intermediate value of hunting cooperation and regain stability via subcritical Hopf-bifurcation if the strength of hunting cooperation is large. 

The existence of multiple limit cycles and bi-stability via subcritical Hopf-bifurcation has been observed by the authors. The solutions tend to a coexisting equilibrium  or oscillate periodically depending on the initial population size. 

When the system shows limit cycle oscillations around the interior equilibrium, increasing the strength of hunting cooperation or cost of fear, then the population oscillations will be replaced by a stable focus. 

It is concluded that fear factor has a stabilizing effect, and by including it, the stable system remains stable and the oscillating system becomes stable by excluding the existence of periodic orbits. It is observed that fear factor has more stabilizing effect than hunting cooperation and makes the system more robust.

For low level of fear, if the cooperation level is increased then the system switches it’s dynamics from stable to limit cycle oscillation and from limit cycle oscillation to stable focus. When the strengths of hunting cooperation and fear factor are very high, prey population goes to extinction due to overexploitation by the predators. In general,   increasing both hunting cooperation and fear factor, the density of both species populations eventually will decrease.

In conclusion, the system with hunting cooperation and fear phenomena exhibits rich dynamical behaviours so these two characteristics are interesting factors with an important role in determining the long term population dynamics.

\color{black}

\vspace{0.5cm}

Insipired by \cite{fearWang}, S. K. Sasmal considered in \cite{Sasmal} a model with a reduction of prey growth rate due to the fear effect and also with the prey species affected by Allee effect. This model is given by
\begin{equation}\label{27}
	\begin{split}
		\dot{x}&=r x \left( 1-\frac{x}{K}\right) (x-m)\frac{1}{1+sy}-axy,\\
		\dot{y}&=a \alpha x y - ny,  
	\end{split}
\end{equation}
where $s$ represents the effect of fear and $m$ the Allee effect.

The authors work on a dimensionless system obtained through variables changes, and obtain the following results for the original system.
The system presents three boundary equilibria which always exist, $E_0=(0,0)$, $E_1=(m,0)$ and $E_2=(K,0)$, and a positive equilibrium $E^*$ which appears under condition 
\begin{equation}
	m < \frac{n}{a \alpha}<s.
\end{equation}
The origin is always locally asymptotically stable. The boundary equilibrium $E_1$ is an unstable node if $m>n/a\alpha$ and it is saddle if $m<n/a\alpha$. The boundary equilibrium $E_2$ is a saddle if $n/a \alpha s<1$ and it is locally asymptotically stable if $n/a\alpha s>1$. The equilibrium $E^*$ is locally asymptotically stable if $(s+m)/2<n/a\alpha<s$ and it is unstable if $m<n/a\alpha<(s+m)/2$.

We note that the origin is  asymptotically stable, so there  exists always some initial conditions that will lead to extinction of both species. Due to strong Allee effects in prey, when the predator invasion is large enough, then predator grows very fast to drive prey below its survival threshold density, and extinction of both the populations is certain.

The analysis carried out concluded that the stability of the positive equilibrium does not depend on the fear effect, the same showed in
\cite{fearWang}, but here the cost of fear changes the density of the predator species on the coexistence equilibrium: as the value of $s$ increases, the density of predators on the positive equilibrium decreases.

The system presents a subcritical Hopf bifurcation at $n/a\alpha=(s+m)/2$, and changing the values of $m$ and $n$ can produce bi-stability or stable oscillatory coexistence of both species. The author observed that modifications on the value of $s$ can only change the density of predator
species at the positive equilibrium, but it does not affect the stability.

\vspace{0.5cm}

Motivated by this work, L.Lai, Z. Zhu and F.Chen, in \cite{fear-alleeaditive} replace the usual Allee effect with the additive Allee effect, which they take from \cite{Alleeaditivo}. They study the system
\begin{equation}\label{aditivo}
	\begin{split}
		\dot{x}&= r x \left( 1- x-\frac{m}{x+b}\right) \frac{1}{1+sy}- a x y,\\
		\dot{y}&= a \alpha x  y - n y, 
	\end{split}
\end{equation}
where the term $1/(1+sy)$ represents the fear effect and the term $m/(x+b)$ represents the additive Allee effect. In particular, the values of $s$ and $m$ represent, respectively, the intensity of fear effect and Allee effect.

Depending on the parameters the system can present two or three boundary equilibrium points. The origin $E_0$ is always an equilibrium for system \eqref{aditivo}. It is a stable node if $b<m$ or $b=m=1$, a saddle-node if $b=m\neq1$ and a saddle is $b>m$. When $b\in(0,1)$ and $b<(b+1)^2/4<m$ there are no more boundary equilibria.

Now we consider the different cases depending on the existence of other boundary equilibria besides $E_0$ and their local stability.
\begin{itemize}
	\item  If $b\in(0,1)$ and  $b<m=(b+1)^2/4$ then  $E_1=\left( (1-b)/2,0\right)$ is an equilibrium for system \eqref{aditivo}. It is a saddle-node if $(\alpha a (1-b) /2)-n\neq0$ and it is a saddle if $(\alpha a (1-b) /2)-n=0$.
	\vspace{0.2cm}
	\item  If $b\in(0,1)$ and  $b<m<(b+1)^2/4$ then there exist two boundary equilibria $E_2=(x_2,0)=\left((1-b-\sqrt{(b+1)^2-4m})/2,0 \right) $ and $E_3=(x_3,0)=\left((1-b+\sqrt{(b+1)^2-4m})/2,0 \right)$ for system \eqref{aditivo}. $E_2$ is always unstable and depending on the sign of $\alpha a x_2 - n$, it can be a saddle, a saddle-node or an unstable node. The equilibrium $E_3$ is a stable node if $\alpha a x_3 - n <0$, a saddle-node if  $\alpha a x_3 - n =0$ and a saddle if  $\alpha a x_3 - n >0$.
	\vspace{0.2cm}
	\item  If $b\in(0,1)$ and $0<m\leq b < (b+1)^2/4$, or $b=1$ and $0<m<1$, or $b>1$ and $0<m<b<(b+1)^2/2$, then $E_3$ is an equilibrium for system \eqref{aditivo} with the same stability as in the previous case. 
\end{itemize}

Besides the boundary equilibria, the system \eqref{aditivo} has a positive equilibrium $E^*=(x^*,y^*)=\left( n/\alpha a, (-b+\sqrt{\Delta})/2as \right)$, with $\Delta=a^2+4asr\left(1-x^*-m/(x^*+b) \right)$, when $1-x^*-m/(x^*+a)>0$.
The equilibrium $E^*$ is unstable if $b<\sqrt{m}-x^*$ and locally asymptotically stable if $b>\sqrt{m}-x^*$, and if also $E_3$ exists and is unstable, the asymptotical stability is global. This allows us to know under what conditions the system always evolves towards the coexistence of both species.

The authors give conditions for the existence of saddle-node bifurcation and transcritical bifurcation of boundary equilibriums. They apply the Sotomayor's theorem, which can be found in \cite{thsotomayor1,thsotomayor2}. The saddle-node bifurcation occurs from $E_1$ at $m=(b+1)^2/4$ if $b<1$ and $\alpha a (1-b)-2n \neq0$. The transcritical bifurcation occurs from $E_0$ at $b=m$ when $m<1$.
They also study the supercritical Hopf bifurcation of the positive equilibrium that occurs when $b=\sqrt{m}-x^*$.

\begin{figure}[H]
	\centering
	\subfloat[Extinction of both species in system \eqref{27} (with $m=2$).]{\includegraphics[width=60mm]{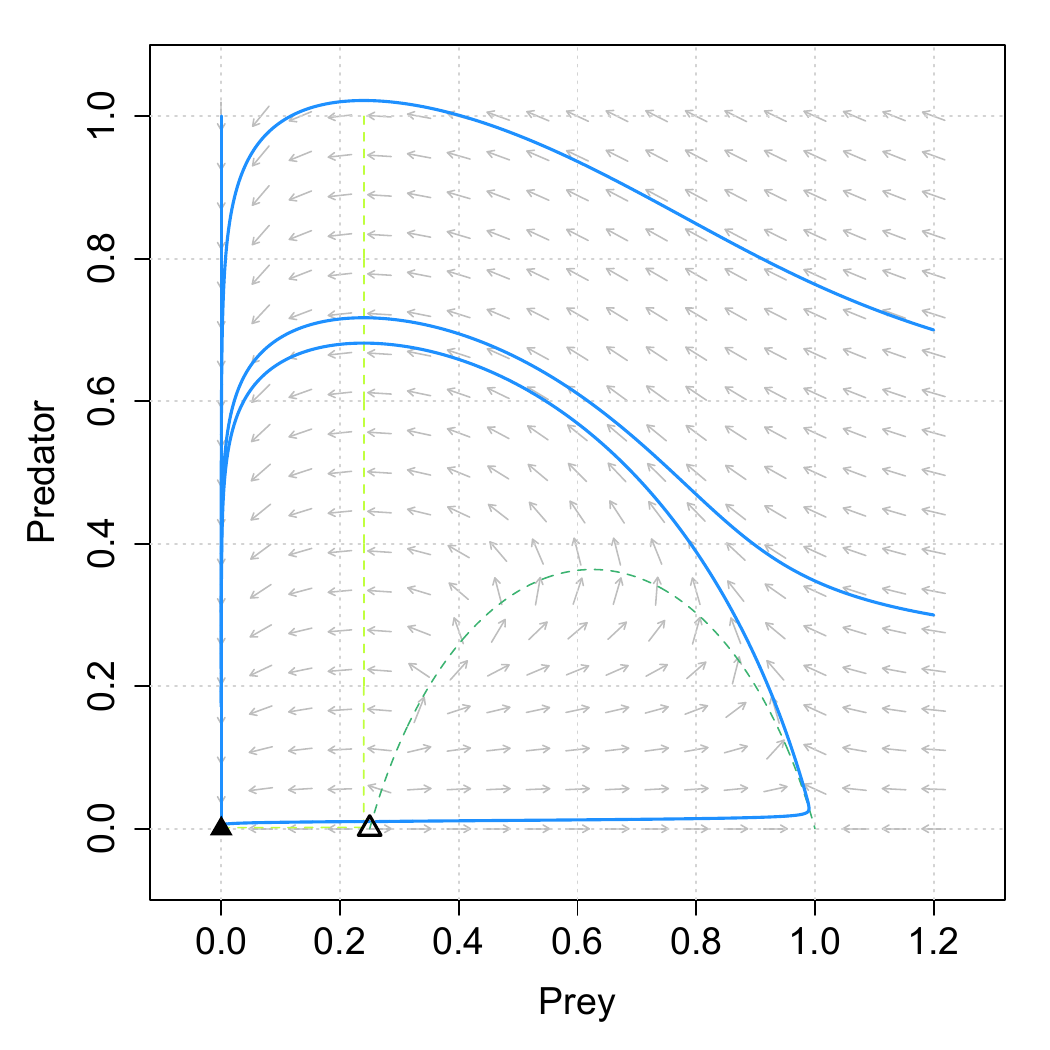}}\\
	\subfloat[Coexistence of both species in system \eqref{aditivo} (with $b=20$).]{\includegraphics[width=60mm]{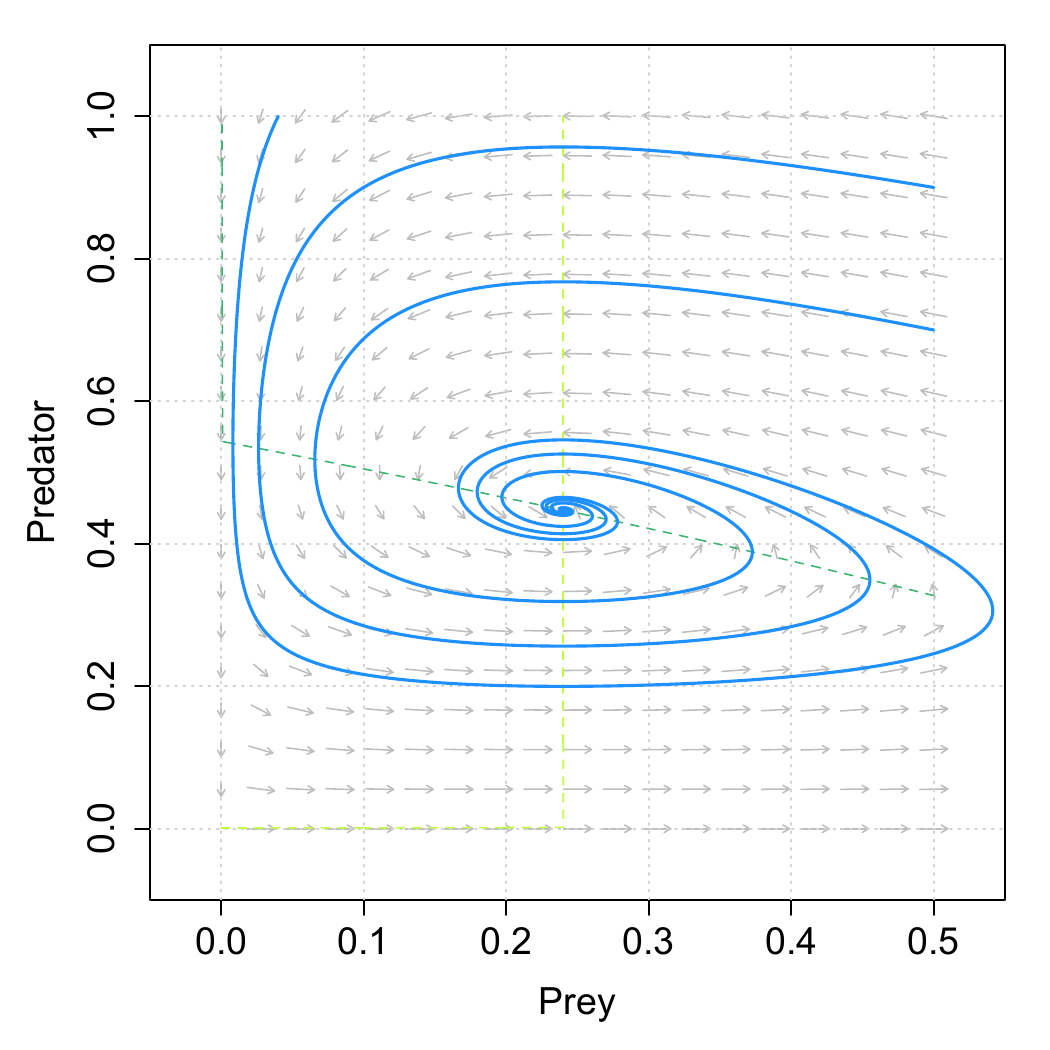}}
	\subfloat[Existence of stable limit cycle in system \eqref{aditivo} (with $b=0.26$).]{\includegraphics[width=60mm]{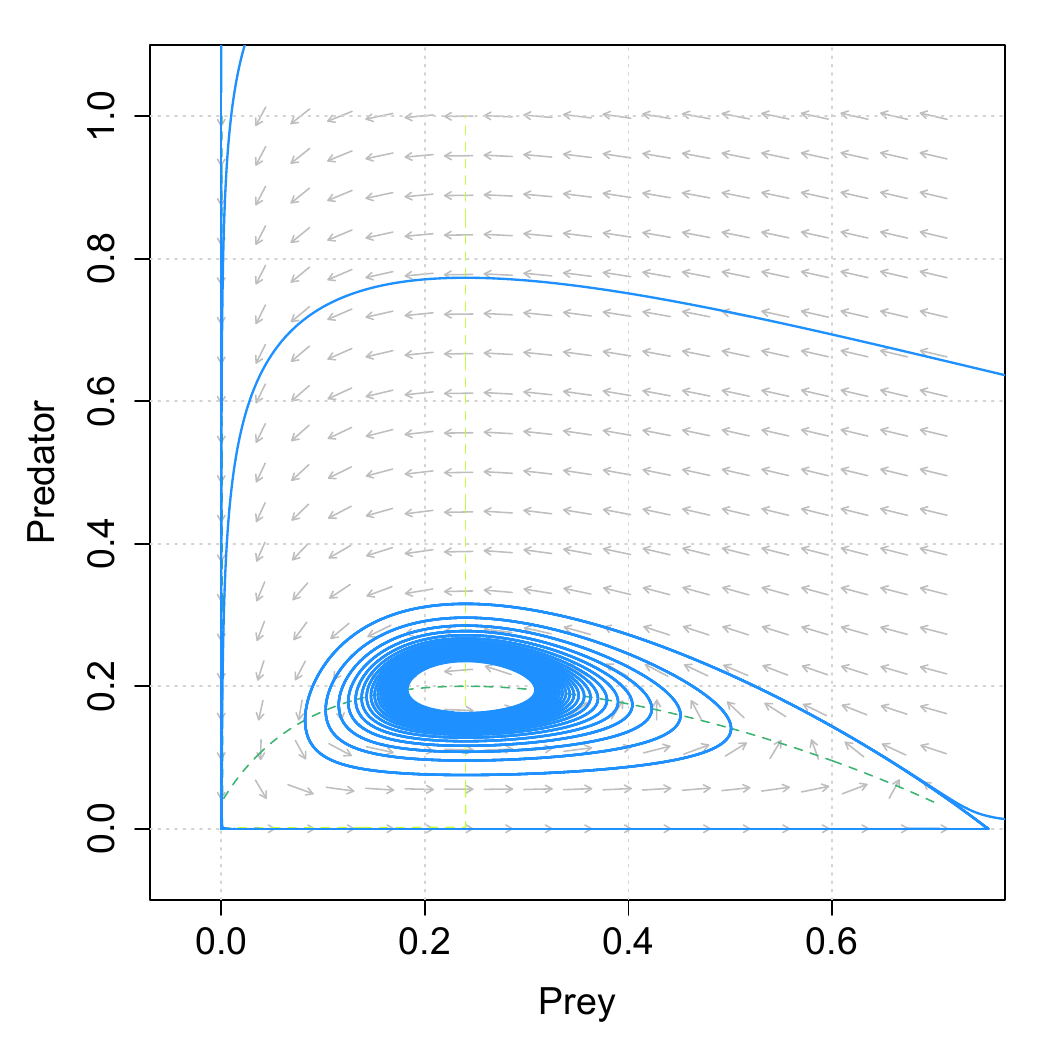}}
	\caption{ Dynamic behaviours with $r=1$, $m=0.25$, $s=1.5$, $a=1$, $K=1$, $\alpha=0.5$ and $n=0.12$.}
	\label{fig:(29)(27)}
\end{figure}

In general, the additive Allee effect gives rise to a more complex dynamics than the multiplicative Allee effect.  As an example, for a fixed set of parameters we show in Figure \ref{fig:(29)(27)} a comparation between the  behaviours of system \eqref{aditivo} and system \eqref{27}. In the system with multiplicative Allee effect, \eqref{27}, both species become extinct while in system with additive Allee effect, \eqref{aditivo}, for the same common parameters and varying the additive Allee effect coefficient $b$ we get different dynamics. 
With regard to the influence of fear effect, this work conclude that it only affects the final density of predator species.

\vspace{0.5cm}

Now we will focus on some models that include a prey refuge.  
First, in  \cite{fear-preyrefuge}, H. Zhang, Y. Cai, S. Fu and W. Wang,  investigated a Holling-II predator prey model incorporating the fear effect and a prey refuge. Their results show that the impact of prey refuge and fear effect is significant. They based their study on the model proposed by T. K Kar in \cite{Kar}, i.e: 
\begin{equation}
	\begin{split}
		\dot{x}&=\alpha x \left( 1- \frac{x}{K}\right)- \frac{\beta(1-m)x y}{1+a(1-m)x},\\[0.1cm]
		\dot{y}&=-\gamma y + \frac{c \beta (1-m)x y}{1+ a(1-m)x},
	\end{split}
\end{equation}
on which they add the fear effect by multiplying by the factor $f(s,y)=1/1+sy$, obtaining the model
\begin{equation}\label{m-pr}
	\begin{split}
		\dot{x}&=\frac{\alpha x}{1 + s y}- b x^2 - \frac{\beta(1-m)x y}{1+a(1-m)x},\\[0.1cm]
		\dot{y}&=-\gamma y + \frac{c \beta (1-m)x y}{1+ a(1-m)x},
	\end{split}
\end{equation}
where $s>0$ refers to the level of fear, and $b=\alpha/K$.
This system represents what they call the anti-predator defense due to fear, and for proposing it, they are also inspired by the results on \cite{fearWang}.

The origin $E_0=(0,0)$ is always an equilibrium for system \eqref{m-pr} and it is always a saddle point. There is also a boundary equilibrium $E_1=(\alpha/b,0)$ that always exists. This equilibrium is locally asymptotically stable if $c\beta-a\gamma \leq 0$, or $c \beta - a \gamma>0$ and $1- b\gamma/\alpha(c\beta-a\gamma)<m<1$, so under these conditions there exists a region from which the system always goes to the extinction of the predator species.
$E_1$ is unstable if $c\beta - a\gamma >0$ and $0\leq m < 1- b\gamma/\alpha (c\beta-a\gamma)$. We note that this conditions do not depend on the parameter $s$, i.e they do not depende on the fear effect. When $E_1$ is unstable, it appears a coexistence equilibrium 
\begin{equation*}
	E^*=(x^*,y^*)=\left( \frac{\gamma}{(c\beta-a\gamma)(1-m)},\frac{1}{2s}\left( \sqrt{\Delta}-1-\frac{kbc\gamma}{(c\beta-a\gamma)^2(1-m)^2}\right) \right),
\end{equation*}
where
$\Delta= \left(sbc\gamma/(c\beta-a\gamma)^2(1-m)^2-1\right)^2+4sc\alpha/(c\beta-a\gamma)(1-m)$.
This positive equilibrium is locally asymptotically stable if 
\begin{equation}\label{H12}
	1-\frac{b(a\gamma-c\beta)}{a\alpha(c\beta-a\gamma)}\leq m < 1 - \frac{b\gamma}{\alpha(c\beta-a\gamma)}  \;\; \text{ or } \;\;
	0 \leq m < 1 - \frac{b(a\gamma+c\beta)}{a\alpha(c\beta-a\gamma)}
\end{equation}
and 
\begin{equation*}
	s > \frac{a(c\beta-a\gamma)^2(1-m)^2(a\alpha(c\beta-a\gamma)(1-m)-b(a\gamma+c\beta))}{ b^2c^2\beta(a\gamma+c\beta)^2}=S^*.
\end{equation*}
The asymptotical stability is global if $0<c\beta-a\alpha\leq a\alpha$ and the first condition in \eqref{H12} holds. At last, $E^*$ is unstable if  the second condition in \eqref{H12} holds and $0<s<S^*$.

The existence of Hopf bifurcation and limit cycles is also studied in this paper. If the second condition in \eqref{H12} holds, then the system has a Hopf bifurcation from $E^*$ at $s=S^*$, and if $0<s<S^*$ then it is proved that there exists a limit cycle. The Poincaré-Bendixon theorem is used in the proof of this result. 

There are some important conclusions, some that agree with other works and others that yield different results. 
When the coexistence equilibrium exists, the increase of the fear parameter $s$ causes the decrease of the predator density, but does not induce the extinction. These facts agree with the results in \cite{fearWang} and \cite{Sasmal}.

Moreover, the fear effect can increase the stability of the system, because it can prevent the existence of limit cycles. When the parameter $m$ related with the prey refuge has a value on the interval $[0,1-b(a\gamma+c\beta)/a\alpha(c\beta-a\gamma))$, the increase of the parameter $s$ changes  the local stability of $E^*$ from unstable to asymptotically stable, and the limit cycle vanishes, so fear makes the system evolve to the coexistence of both species. This result does not agree with the ones in \cite{Sasmal}.

With regard to the prey refuge, when $m\in[0,1-2b\gamma/\alpha(c\beta-a\gamma))$, with the increase of $s$ the predator density changes from increase to decrease. When $m$ reaches a high risk threshold of the prey refuge, the predator goes to extinct. 
\vspace{0.5cm}

We have seen how incorporating different biological issues, such as the Allee effect and the prey refuge, can change the qualitative behaviours. Now we continue combining some of these biological characteristics, in particular, the next model incorporates at the same time the Allee effect and the prey refuge.
Motivated by the work of Sasmal, in \cite{pr-fear-allee} the authors analyze the influence of prey refuge on the model raised in \cite{Sasmal}. They study the system
\begin{equation}\label{33}
	\begin{split}
		\dot{x}&= r x \left( 1 - \frac{x}{K}\right) (x- m)\frac{1}{1+s y}- a (1-\eta)xy,\\
		\dot{y}&= a \alpha (1-\eta)xy - n y,
	\end{split}
\end{equation}
where $\eta\in (0,1)$ is the prey refuge constant and so $\eta x(t)$ is the capacity of the refuge at time $t$. 
If we compare with the previous model included in this review, we are adding the Allee effect in addition to fear effect and prey refuge.

For system \eqref{33} the origin is always a equilibrium, particularly, a stable node.  There are also two boundary equilibria, $E_1=(K,0)$ and $E_2=(m,0)$. $E_1$ can be a saddle, an attracting saddle-node or locally asymptotically stable, depending on the relation between the parameter $\eta$ and $n/a\alpha b$. In a similar way, $E_2$ can be a repelling saddle-node or a saddle depending on the relation between the parameter $\eta$ and $n/a\alpha m$. 

Under certain conditions there exist also a positive equilibrium $E^*$. 
As they consider the prey refuge as a parameter, they arrive to a threshold condition for the stability of the system. It is also proved that the system has a supercritical Hopf bifurcation. 

\begin{figure}[H]
	\centering
	\subfloat[Possible coexistence in system \eqref{27}.]{\includegraphics[width=55mm]{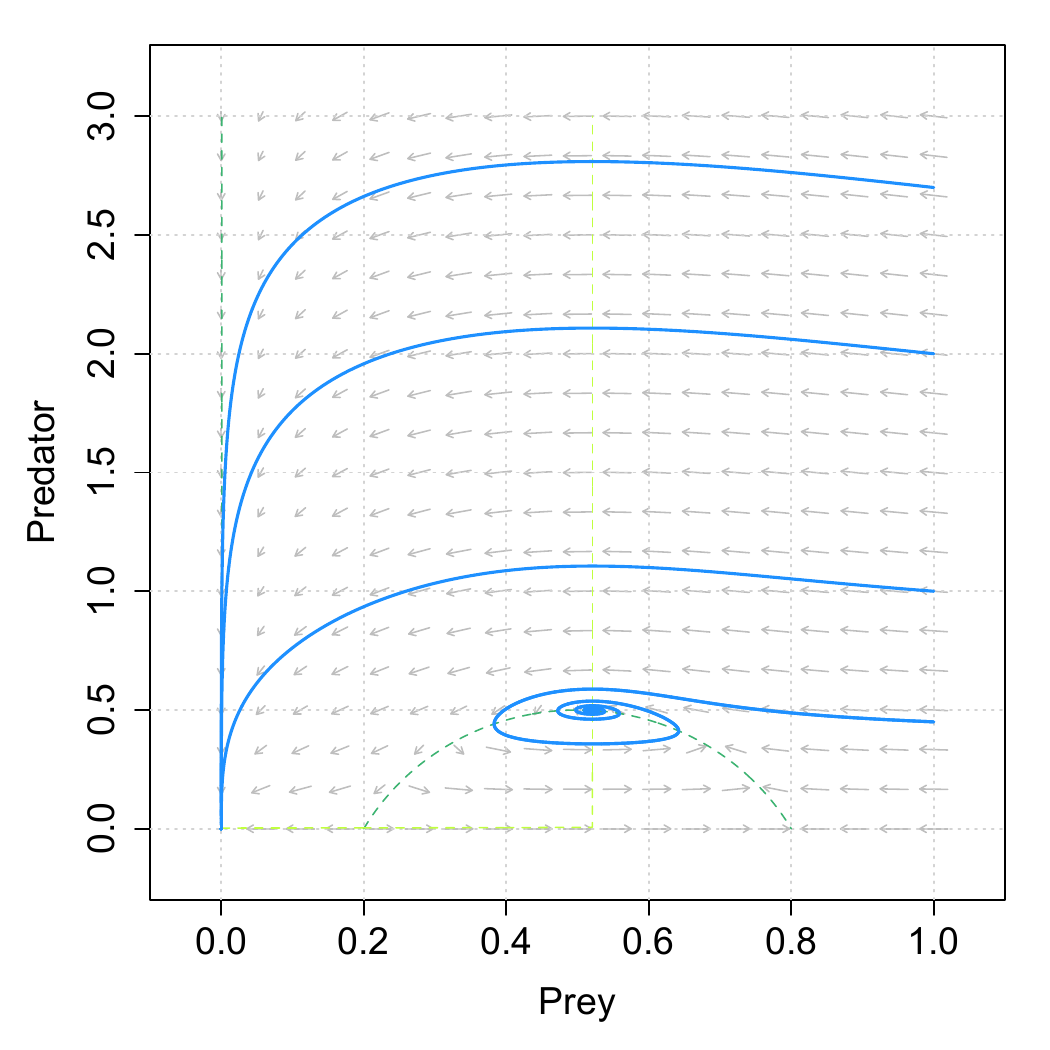}}
	\subfloat[Extinction of predators in system \eqref{33}.]{\includegraphics[width=55mm]{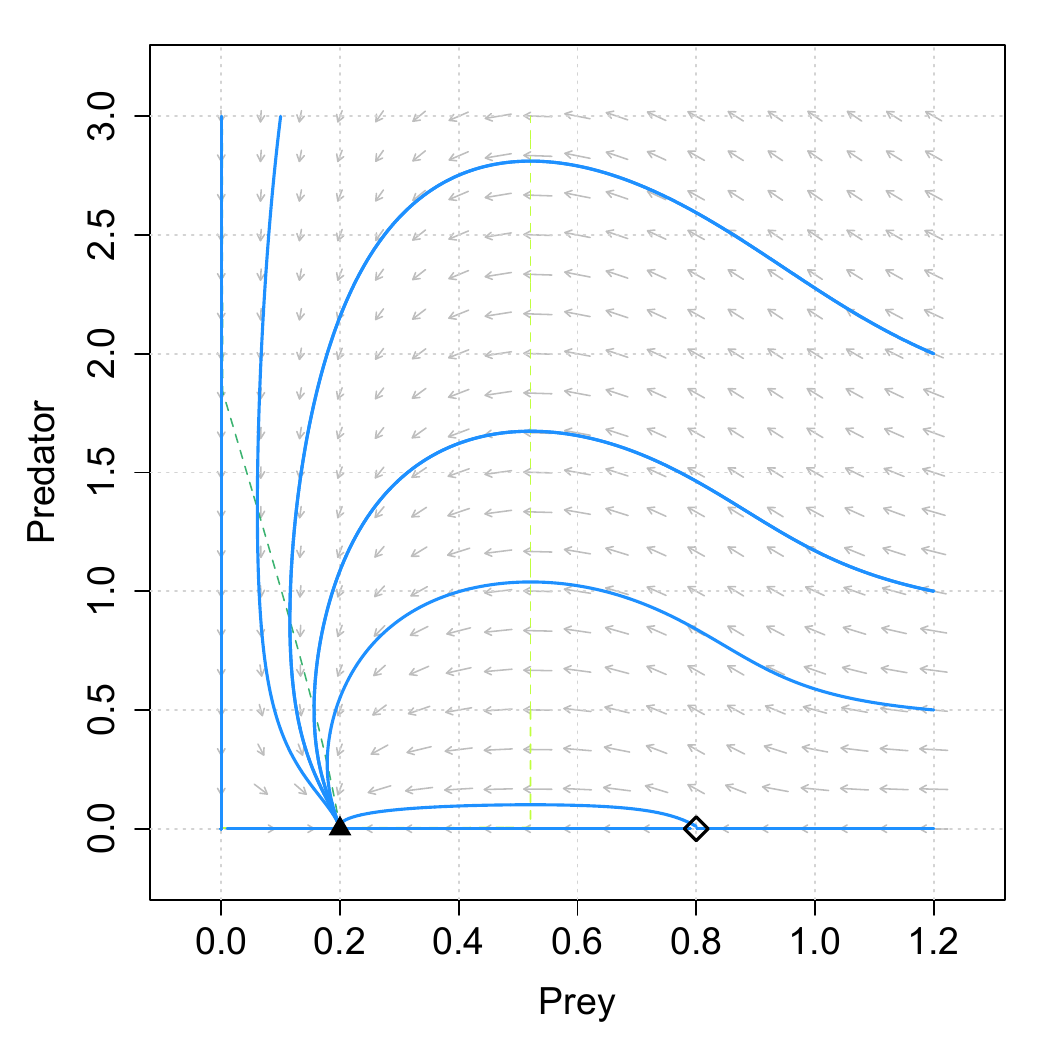}}
	\caption{Dynamic behaviours with $r=2$, $K=0.2$, $m=0.8$, $s=1$, $a=1.2$, $\alpha=0.8$, $n=0.5$ and $\eta=0.8$.}
	\label{fig:(33)(27)}
\end{figure}

It is concluded that the prey refuge has an important role on the dynamics of the system. Comparing the results with the ones obtained in \cite{Sasmal}, the dynamic is more complex in this case. The influence of the prey refugee may not be positive for the population, as shown in Figure \eqref{fig:(33)(27)}. For a fixed set of parameters, in system \eqref{27} without prey refugee the coexistence of both species is possible, but in system \eqref{33} it is not possible.
The fear effect and Allee effect do not affect the final density of prey on the coexistence equilibrium, but they can cause a decrease in the predators population.

\section{Cannibalism}\label{sec:cannibalism}

Cannibalism is the act of killing and at least partial consumption of conspecifics. This actively occurs in more than 1300 species in the nature \cite{1300can}. In the literature, cannibalism has been principally considered in the predator species, and in general the results agree that this cannibalism stabilizes the system and causes persistence in a population doomed to go extinct. The first contribution in this sense was the work of Kohlmeir \cite{Kohlmeier}, who considered cannibalism in the predator in the classic Rosenzweig-McArthur model. As an example we present the work of H. Deng, F. Chen, Z. Zhu and Z. Li \cite{cannibalism} and after we present the first work that considers the effect of prey cannibalism in a two species unstructured ODE model \cite{cannibalism2}.

\vspace{0.5cm} 
In \cite{cannibalism}, the authors consider a modification of the original Lotka-Volterra model which includes the cannibalism in the predator species. The system they propose is the following
\begin{equation}\label{34}
	\begin{split}
		\dot{x}&=x(b-\alpha x - my), \\
		\dot{y}&=y (-\beta+c_1+nx)-\frac{cy^2}{y+d},
	\end{split}
\end{equation}
where $c_1<c$ and $cy^2/(y+d)$ denotes the cannibalism of the predator.

The authors obtain the equilibria of the system and study their local and global stability. The results are compared with the ones obtained for the original system without cannibalism.
There are two equilibria that always exist, the origin $E_0$ and the boundary equilibrium $E_1=(b/\alpha,0)$, and another boundary equilibrium $E_2=(d(c_1-\beta)/(\beta+c-c_1)$ that exists if $c_1>\beta$. A coexistence equilibrium $E^*=(x_2^*,y_2^*)$,  exists if $0<\beta-c_1<bn/\alpha$ or if $\beta \leq c_1 < \beta + bc/(b+dm)$, where
\begin{equation*}
	\begin{split}
		&x_2^*= \frac{B-\sqrt{B^2-4AC}}{2A}, \hspace{1.5cm}  y_2^*=\frac{b-\alpha x_2^*}{m}, \\[0.1cm]
		A= \alpha n, \;\;  \;  &B=\alpha(\beta+c-c_1)+bn+dmn, \; \text{ and }  \; C=b(\beta+c-c_1)+dm(\beta-c_1).
	\end{split}
\end{equation*}
This positive equilibrium $E^*$ is globally asymptotically stable whenever it exists.

The equilibrium $E_0$ is a saddle if $c_1<\beta$, a saddle-node if $c_1=\beta$ and an unstable node if $c_1>\beta$. The equilibrium $E_1$ is a saddle if $c_1\geq\beta$ or if  $0<\beta-c_1<bn/\alpha$; it is a stable node if $\beta-c_1>bn/\alpha$ and a saddle node if $\beta-c_1=bn/\alpha$. At last, when $E_2$ exists, it is a saddle if $\beta<c_1<\beta + (b(\beta+c)+\beta dm)/(b+dm)$,
a stable node if  $c_1>\beta + (b(\beta+c)+\beta dm)/(b+dm)$, 
and a saddle node if $c_1=\beta + (b(\beta+c)+\beta dm)/(b+dm)$.

The conclusion is that the presence of cannibalism modifies the stability of the system, but not always in the same direction.
When in the original system the boundary equilibrium in which there are only prey is globally asymptotically stable, the introduction of cannibalism (in a suitable quantity) benefits the system, i.e, it allows the two species to coexist without  becoming extinct, since a globally asymptotically stable positive equilibrium appears.
If the parameter related with cannibalism increases too much, this positive equilibrium disappears and another boundary equilibrium appears, with the result that  prey becomes extinct but even so the predators manage to survive on cannibalism.

On the other hand, when the original Lotka-Volterra system is in conditions in which the coexistence equilibrium exists and is globally asymptotically stable, as the parameter associated with cannibalism increases the population of prey increases and the population of predators decreases. In this way,  when a certain value of this parameter is reached, cannibalism causes the extinction of the predators. 

\vspace{0.5cm}

In \cite{cannibalism2}, A. Basheer, E. Quansah, S. Bhowmick and R. D. Parshad consider the 
the first ODE model with two species and prey cannibalism. They compare the simple Holling-Tanner model with ratio dependent functional response in prey \cite{HollingTanner}, i.e, the system
\begin{equation}\label{HT}
	\begin{split}
		\dot{x}&=x(1-x)-\frac{xy}{x+\alpha y}, \\
		\dot{y}&=\delta y \left(\beta - \frac{y}{x} \right) , 
	\end{split}
\end{equation}
with the systems obtained by adding cannibalism in the predator or in the prey species. For the original system, if $\beta<1/(1-\alpha)$ then there exists a positive equilibrium  $(x^*,y^*)=(1-\beta/(1+\alpha \beta), \beta(1-\beta/(1+\alpha\beta)))$, and it is stable whenever it exists.

By considering the inclusion of predator cannibalism they obtain the system 
\begin{equation}\label{HTpredator}
	\begin{split}
		\dot{x}&=x(1-x)-\frac{xy}{x+\alpha y}, \\
		\dot{y}&=\delta y \left(\beta_1 - \frac{y}{\gamma x + c y} \right). 
	\end{split}
\end{equation}
Here the food source of predators is represented by $\gamma x + cy$, where $cy$ is the cannibalism term. There is a gain of energy to the cannibalistic predators from the act of cannibalism. This results in an increase in reproduction which is represented by adding a $\beta_1 y$ term in the predation equation, assuming $\beta_1>\beta$.  For system \eqref{HTpredator} there is a positive equilibrium 
$(x^*,y^*)=(1-\beta_1\gamma/(\alpha \beta_1 \gamma - c \beta_1 + 1), \gamma \beta_1/(1-\beta_1 c))$ if $\gamma<1$, $1-\beta_1 c>0$ and $\alpha\gamma-\gamma + 1/\beta_1>c$. According with the results obtained in \cite{cannibalism}, for this equilibrium the following result is proved: For a parameter set such that the positive equilibrium of system \eqref{HT} is unstable, if $\gamma<\beta/\beta_1$ then there exists a value of the cannibalism rate $c$ for which the positive equilibrium of system \eqref{HTpredator} is stable. This implies that in that case the equilibrium can be stabilised via predator cannibalism. This does not occur in the case $\gamma>\beta/\beta_1$. 

Now, with the inclusion of prey cannibalism instead of predator cannibalism, the authors study the system 
\begin{equation}\label{HTprey}
	\begin{split}
		\dot{x}&=x(1+c_1-x)-\frac{xy}{x+\alpha y}-c\frac{x^2}{x+d}, \\
		\dot{y}&=\delta y \left(\beta - \frac{y}{x} \right). 
	\end{split}
\end{equation}
The cannibalism term $C(x)=c x (x/(x+d)))$ is added in the prey equation. The gain of energy to the cannibalist prey and the consequent increase in reproduction is modeled by adding a $c_1x$ term to the prey equation. It is assumed that $c_1<c$ as it takes depredation of a number of prey by the cannibal to produce one new offspring. For system \eqref{HTprey} the positive equilibrium is
\begin{equation*}
	(x^*,y^*)=\left( \frac{-(m+d+c-c_1-1)+\sqrt{(m+d+c-c_1-1)^2+4d(1-m)+4dc_1}}{2}, \beta x^* \right) 
\end{equation*}
with $m=\beta/(\alpha \beta + 1)<1$ and under some of the following feasibility conditions: 
\begin{enumerate}
	\item $m<1+c_1$ and $m+d+c\geq 1+c_1$, 
	\item $m<1+c_1$ and $m+d+c< 1+c_1$.
\end{enumerate}
The authors answer two questions about the effect of prey cannibalism. Firstly, can prey cannibalism destabilize the positive equilibrium? They prove that there always exists a cannibalism rate $c$ such that the positive equilibrium of system \eqref{HT} is stable while the positive equilibrium of the system \eqref{HTprey} is unstable, so it has a destabilizing effect.
The following question is, can prey cannibalism stabilize the positive equilibrium? In this case it is proved that, if $m+d+x<1+c_1$, then no amount of prey cannibalism can stabilize the unstable positive equilibrium of the system \eqref{HT}.  This corresponds with the first feasibility condition, under which the result is rigorously proven. The authors conjecture the same result for the second feasibility conditions, based on their numerical study. 
Also by numerical analysis the authors note that even with stochasticity in the cannibalism parameter, prey cannibalism is unable to stabilise the dynamics. 

Moreover, the authors study the existence of limit cycles. In particular, they show that, for system \eqref{HTprey}, if the conditions
\begin{equation*}
	\alpha \beta + 1 > \beta  \; \text{ and } \; \left( \frac{\beta}{(1+\alpha\beta)^2}-x^*-\beta \delta\right) \frac{(x^*+d)^2}{dx^*}  >c. \color{black}
\end{equation*}
hold, then there exists at least one limit cycle. Note that the uniqueness was not proved. 

From an ecological point of view, this results show that in a predator-prey community with cyclical oscillations, introducing cannibalism in the prey would just maintain the behaviour, and the population cannot be driven to a stable steady state.

\vspace{0.5cm}

Cannibalism on prey had previously been considered in multicomponent structured model via integer differential equations \cite{canp-ej1}, structured discrete model \cite{canp-ej2} or three species structured ODE models \cite{canp-ej3}.

Also, on mathematical literature about cannibalism, another types of predator-prey models with inclusion of cannibalism have been investigated, as three species ODE models incorporationg stage structure, two species PDE models, discrete models or integro-differential equations models. \cite{can-ej1,can-ej2,can-ej3,can-ej4,Kohlmeier}.

\section{Models with immigration}\label{sec:immigration}

Most predator-prey systems in the wild are not isolated, so it is important to consider the effects of the presence of some number of immigrants. We will focus on studying the inclusion of immigration in the Rosenzweig-MacArthur model and in the Lotka-Volterra model. 

In the literature one can find works of a very different nature that consider immigration. On many occasions, species of prey and predators migrate instantaneously after facing a hostile situation. However, after going through one of these situations, the individuals may encounter some barrier that causes delayed migration. In this sense, some models based on delay equations can be found, \cite{inm_delay,delay2,delay3,delay4}. Some other studies consider the species distributed between different patches, so systems in $\mathbb{R}^n$ with $n>2$ are studied \cite{patches1,patches2,patches3,patches4}.

\vspace{0.5cm}

Taking into account the hypothesis that spatial interaction contributes to the regulation of predator-prey populations, in \cite{inm1} the authors constructed 
a model for spatially interacting populations, adding constant immigration of prey to the Rosenzweig-MacArthur predator-prey model with a Holling type II functional response. They consider the system

\begin{equation}\label{inm1}
	\begin{split}
		\dot{x}&= r x \left( 1- \frac{x}{K}\right) - \frac{x y}{a+x} + c, \color{black}\\
		\dot{y}&= y \left( \frac{\mu x }{a+x}- d\right),
	\end{split}
\end{equation}
where the prey inmigration is given by $c\geq0$.

In the case $c=0$, i. e, for the Rosenzweig-MacArthur model, the dynamics are well known. The origin is always an equilibrium and also it is the point $(K,0)$. There exists a positive stable equilibrium when $\mu>d$ and $a(\mu+d)/(\mu-d)>K>ad/(\mu-d)$, which bifurcates at $K=a(\mu+d)/(\mu-d)$. When $K>a(\mu+d)/(\mu-d)$ a limit cycle emerges. When the carrying capacity $K$ increases, the cycle brings one of both populations closer to zero. 

In the case with immigration, i.e when $c\neq0$, the origin is not an equilibrium point and there exists a boundary equilibrium $E_1=(K/2+\sqrt{(K/2+cK/r))})$. 
Under conditions $\mu>d$ and $K>r\lambda^2/(r\lambda+c)$ there exists a positive equilibrium $E^*(\lambda,\mu (r\lambda (1-\lambda/K)+c)/ d)$. Under any other conditions there is not positive equilibrium and neither limit cycles. 

Although $E^*$ exists, the system does not always have limit cycles, in fact there are no limit cycles if $\mu>d$ and $ r\lambda^2/(r\lambda+c)< K \leq \lambda$. There exists a unique limit cycle if $\mu>d$, $K>\lambda$ and 
\begin{equation}\label{bif}
	(\mu+d)\lambda + \frac{bK(\mu-d)}{r\lambda}<dK.
\end{equation}
When $E^*$ exists, it is globally asymptotically stable if and only if condition \eqref{bif} does not hold.

This allows the authors to state some conclusions from the point of view of ecology.
They assume the hypothesis $\mu>d$, which they consider reasonable for the maintenance of the predator population.
If $K\leq \lambda$ the model without immigration has no positive equilibrium, and the boundary equilibrium is asymptotically stable, so the predators become extinct. 
By incorporating immigration into the model and increasing the value of $c$, a positive equilibrium appears when $K = r \alpha^2/ (r\lambda +c)$. This equilibrium is asymptotically stable, so authors conclude that immigration allows the coexistence.

If $K>\lambda$ there exist a positive equilibrium (for any value of $c$).  Depending on  the values of the parameters, the system can present a  globally aymptotically stable equilibrium, in which case the two population coexist, or it can have a limit cycle, in which case the two populations have an oscillatory behaviour.

For system \eqref{inm1} the authors give necessary and sufficient conditions for the asymptotic stability of the positive equilibrium and the uniqueness  of limit cycles. They use the Liénard's theory, after transforming the system into one of Liénard type by means of a standard way, which is detailed in the proof of the result. Their results clarify that the prey constant immigration dampens the large amplitude of oscillations emerging around the positive equilibrium. This agrees with the general idea that immigration stabilizes predator-prey populations.

\vspace{0.5cm}

In \cite{inm2} the authors consider the modified system with few predator and prey immigration. They analyze the asymptotic stability of the predator-prey system but adding an inmigration factor $c_1(x)$ into the prey population and $c_2(x)$ into the predator population, so they study the system 
\begin{equation}\label{40}
	\begin{split}
		\dot{x}= r x - a x y + c_1(x), \\
		\dot{y}= b x y - n y + c_2(y), 
	\end{split}
\end{equation}
where the immigration functions $c_i$ can be modeled into two ways, as $c_i(z)=c_i$ or as
$c_i(z)=c_i/z$ for  $z\neq0$.
The first expression represents a constant number of immigrants that arrive into a population, which can occur if the habitat is attractive for its quality. If we choose the second expression the number of individuals arriving depends on the population, for high populations densities the number of immigrants is lower.
We recall that the case with $c_i<0$ can be also considered.

Three different models are considered, which can be written in general as
\begin{equation}\label{41}
	\begin{split}
		\dot{x}=r x - \frac{a x^{1+\alpha}y}{1+h x^{1+\alpha}}+c_1(x),\\[0.1cm]
		\dot{y}=\frac{b x^{1+\alpha} y}{1+hx^{1+\alpha}}-ny + c_2(y),
	\end{split}
\end{equation}
where $h$ and $\alpha$ are the functional response coefficients. Note that for $h=0$ and $\alpha=0$ we have a type I functional response, for $h\neq0$ and $\alpha=0$ we have type II functional response and for $h\neq0$ and $\alpha>0$ a type III functional response. 

For each one of these three models, four different types on small immigration are considered, determined by the definition of $c_1$ and $c_2$. Respectively they consider: $c_1(x)=c_1$ and $c_2(y)=0$; $c_1(x)=c_1/x$ and $c_2(y)=0$; $c_1(x)=0$ and $c_2(y)=c_2$; and at last, $c_1(x)=0$ and $c_2(y)=c_2/y$.

For the model with type I functional response, i.e, the linear model,  the authors find that without immigrants the positive equilibrium is unstable and there exist limit cycles. In this case, the coexistence equilibrium becomes stable with any of the four immigration types.

For the model with  type II functional response, the positive equilibrium is unstable without immigrants, but there are no limit cycles. In this case, the coexistence equilibrium becomes locally asymptotically stable with any of the four immigration types.

Finally, for the model with type III functional response, the positive equilibrium is locally asymptotically stable if $rb-hnr>0$. This equilibrium remains locally asymptotically stable regardless of the type of immigration considered.

\begin{figure}[H]
	\centering
	\subfloat[Dynamic behaviour with no immigration.]{\includegraphics[width=50mm]{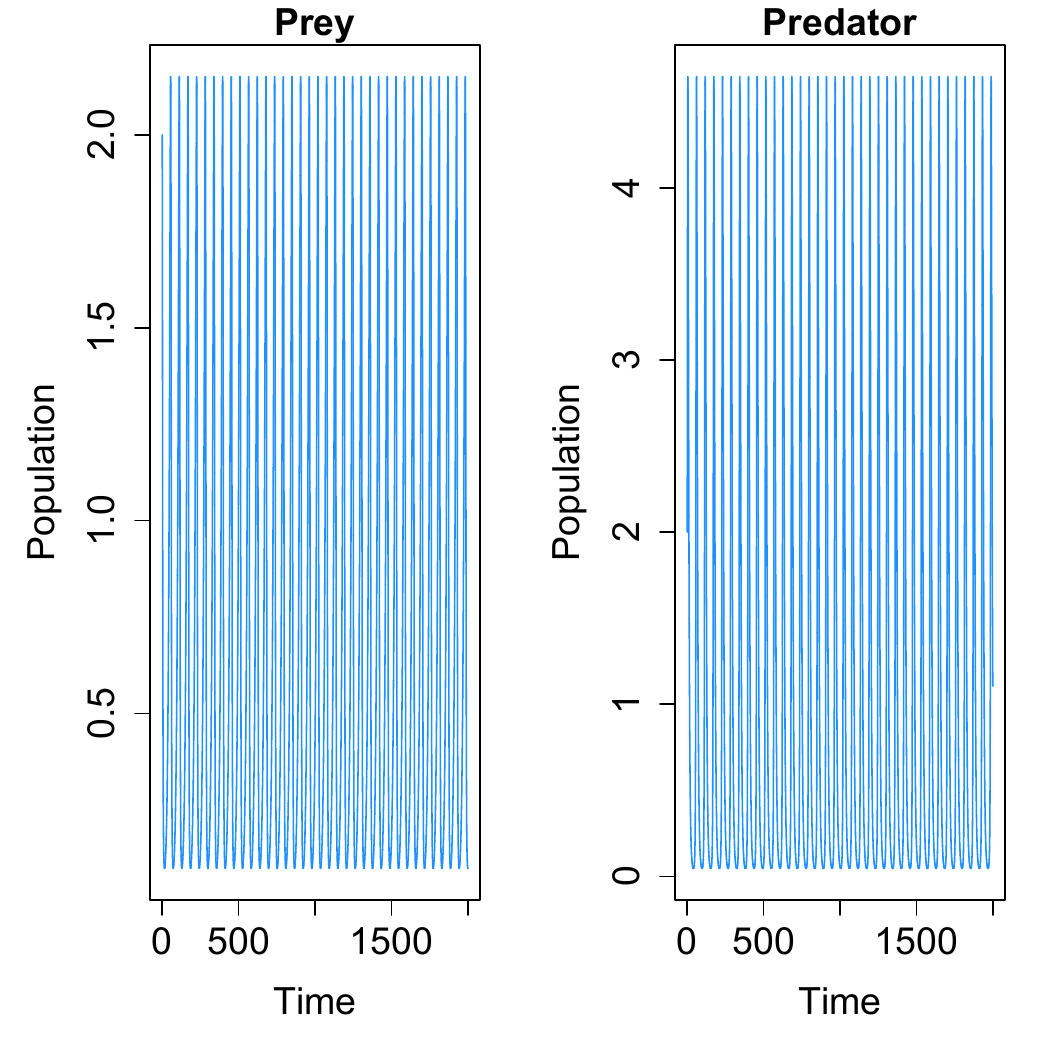}\includegraphics[width=50mm]{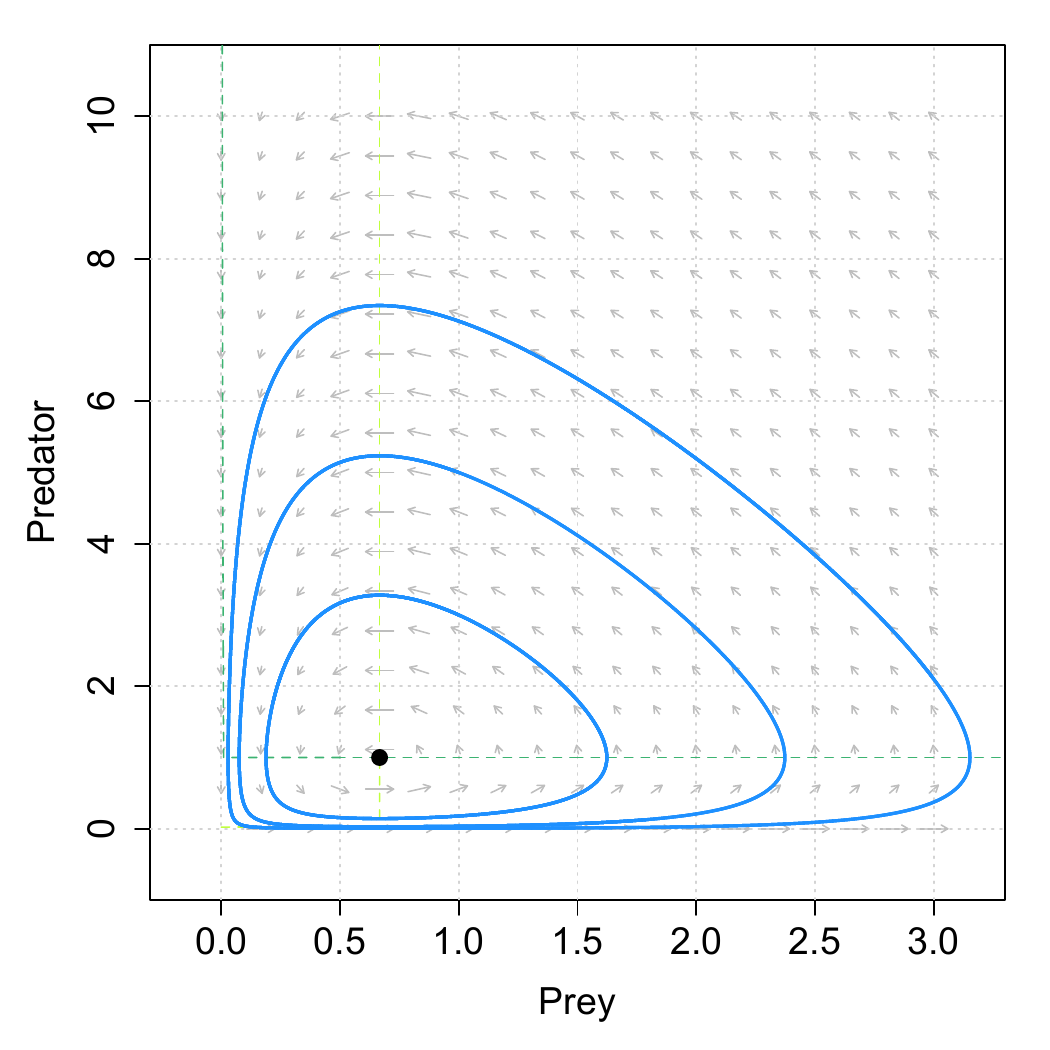}}\\
	\subfloat[Dynamic behaviour with prey immigration ($c_1(x)=0.01$, $c_2(y)=0$).]{\includegraphics[width=50mm]{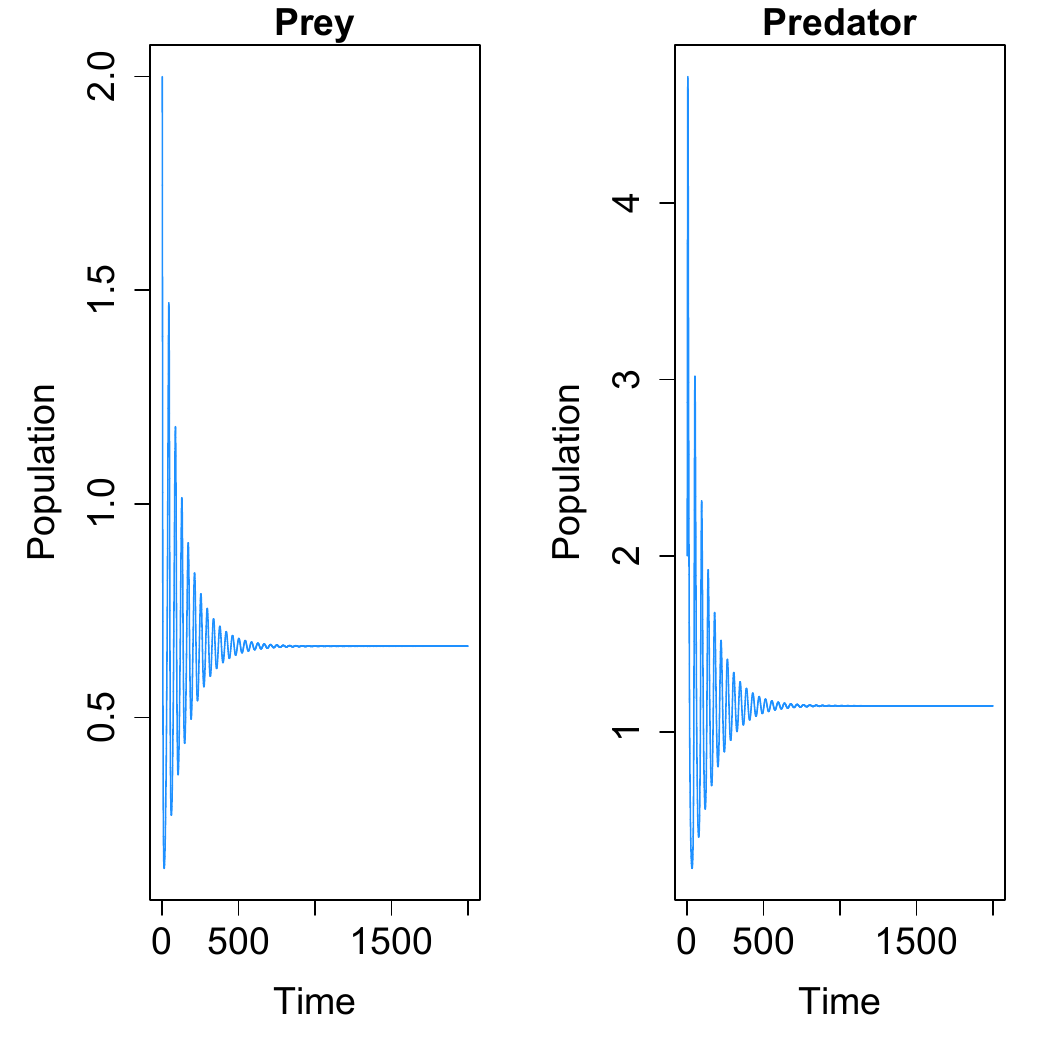}\includegraphics[width=50mm]{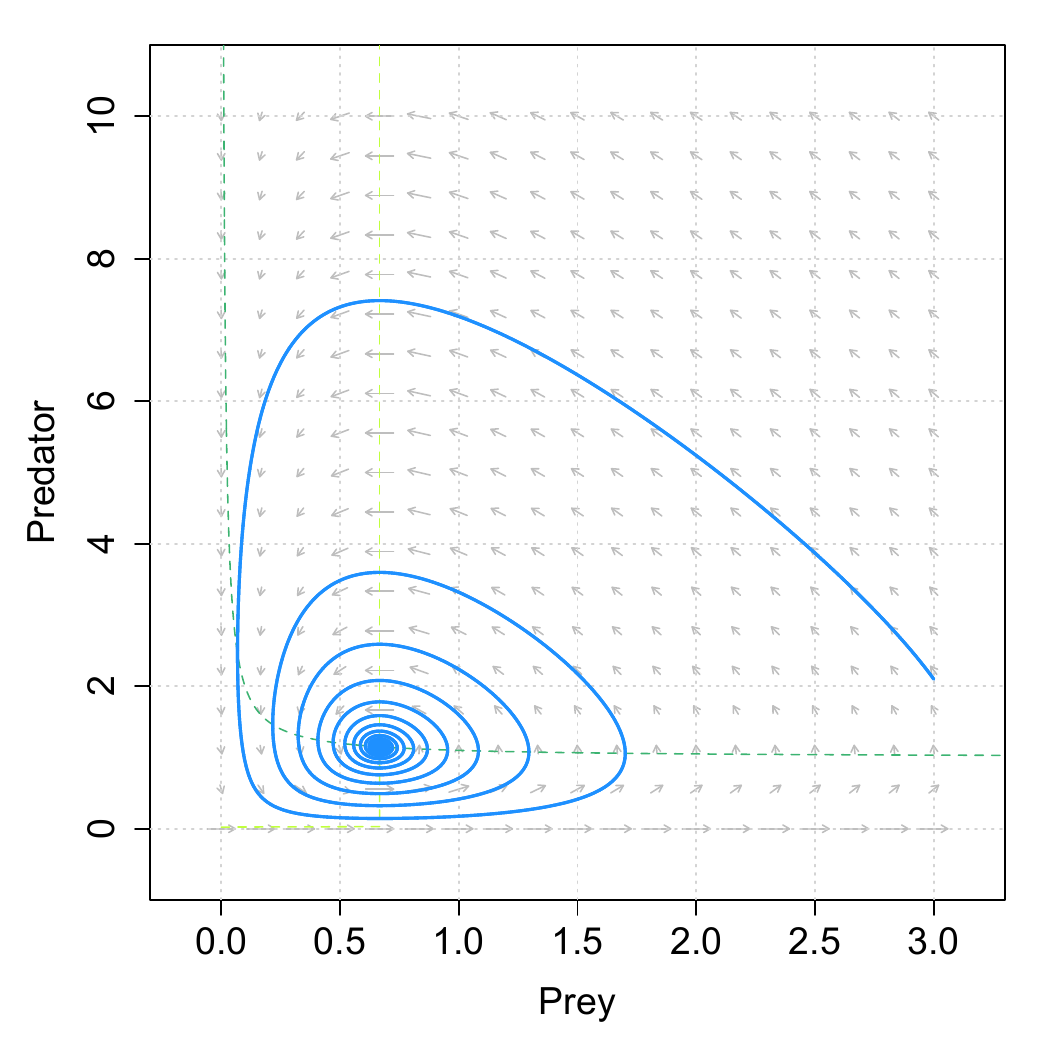}}\\
	\subfloat[Dynamic behaviour with predator immigration ($c_1(x)=0$, $c_2(y)=0.01$).]{\includegraphics[width=50mm]{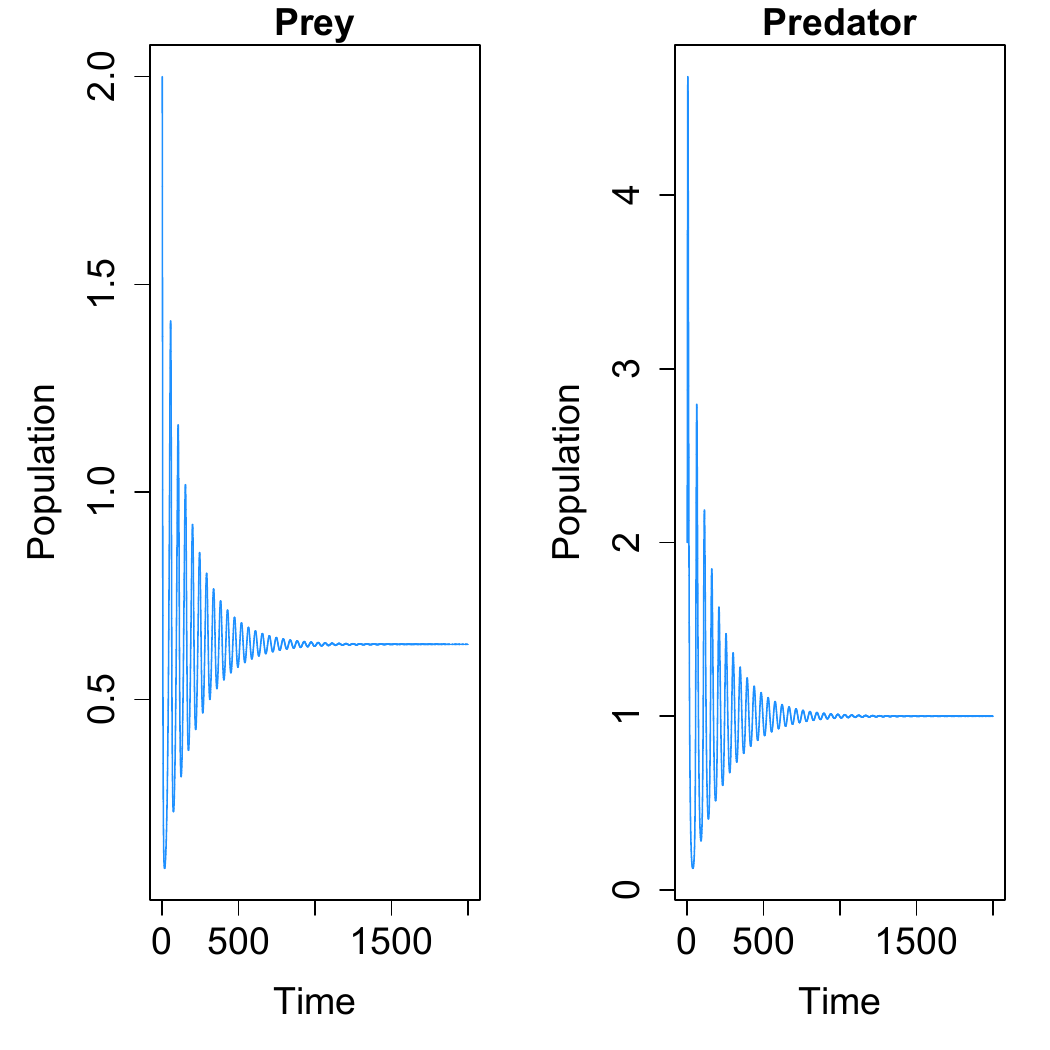}\includegraphics[width=50mm]{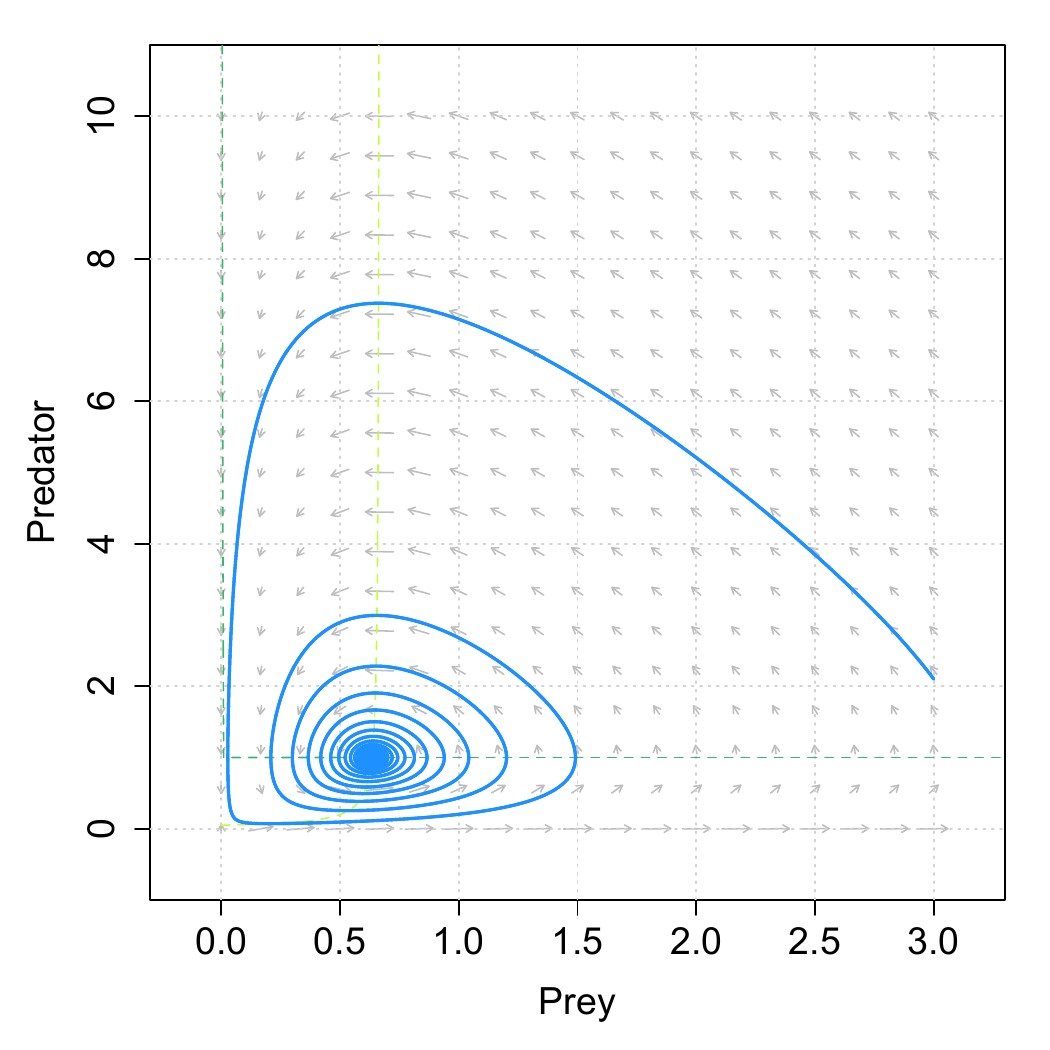}}\\
	\subfloat[Dynamic behaviour with prey and predator immigration ($c_1(x)=0.01$, $c_2(y)=0.01$).]{\includegraphics[width=50mm]{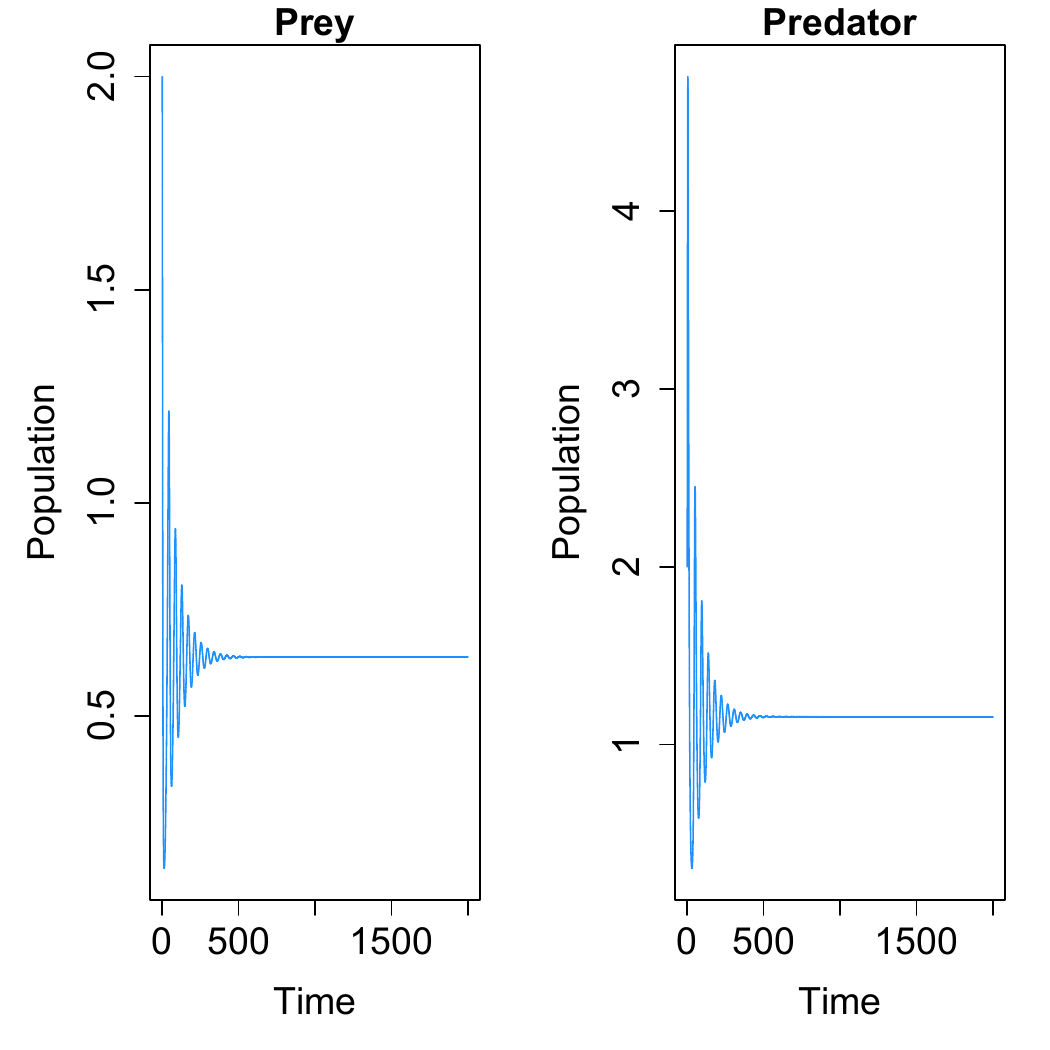}\includegraphics[width=50
		mm]{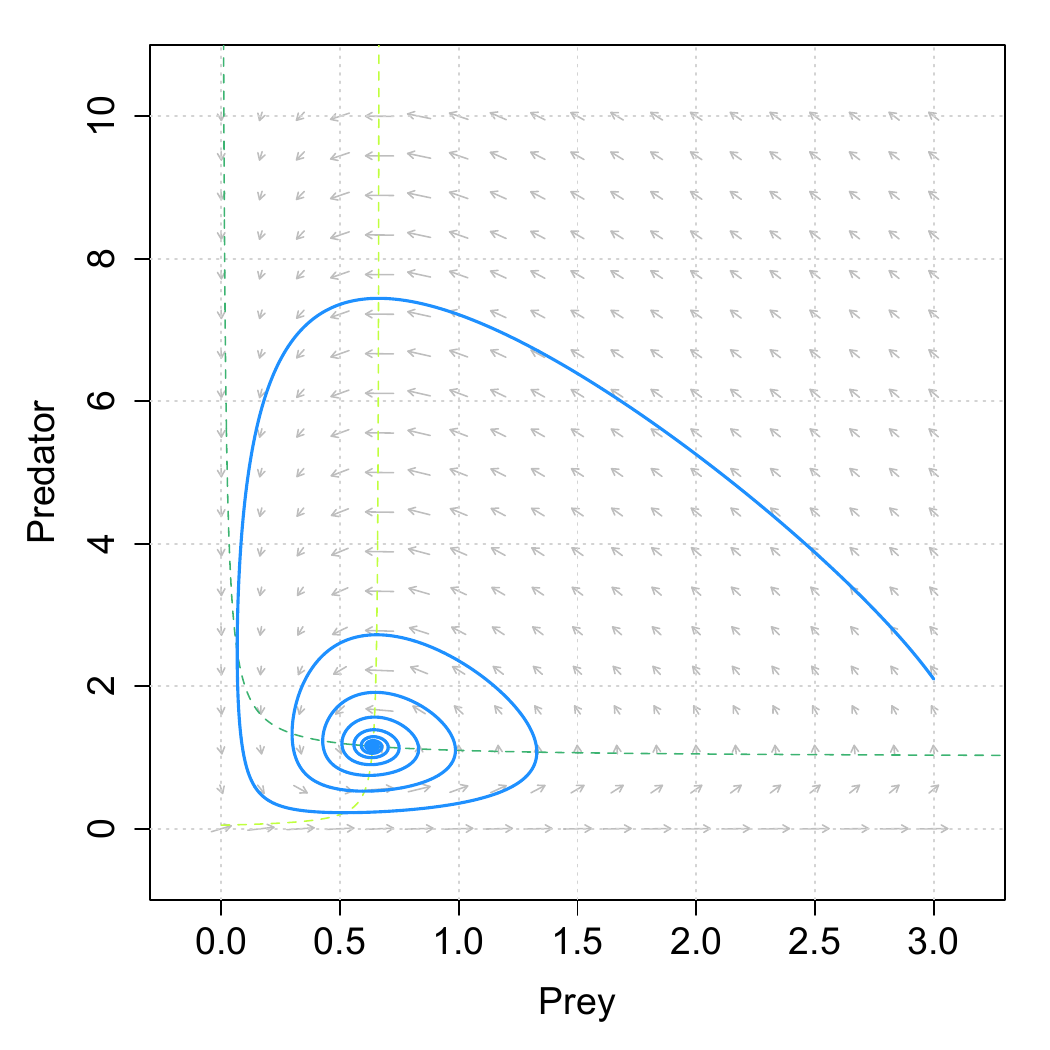}}\\
	\caption{Solutions and phase portraits of system \eqref{40} with parameters $r = 0.1$, $a=0.1$, $b=0.3$ and $n=0.2$.}
	\label{fig:(40)}
\end{figure}

This work concludes that very small immigration into either prey or predator population acts as a stabilizing factor to the Lotka-Volterra system. Adding positive immigration will average out all fluctuations in both the populations of prey and predator. Note also that a positive immigration factor is enough to change the qualitative dynamics of the population. This paper may imply that cyclic populations can be stabilized by adding few immigrations into them. As an example, for a fixed set of parameters, in Figure \ref{fig:(40)} we compare the behaviour of system \eqref{40} with and without constant immigration. 
In most natural populations there are at least a few immigrants over time, and these small number of immigrating individuals are sufficient for obtaining asymptotic stability in the predator-prey system.  Note that the case with immigration in both species is not studied in \cite{inm2}, but we have carried out some simulations and the results obtained are similar to the obtained for the cases with immigration in only one species, as can be seen in Figure \ref{fig:(40)}(d). \color{black}

\section{Conclusions}\label{sec:conclusions}

Population models, and in particular predator-prey models, are a widely studied topic in biomathematics. We have reviewed some recent works focusing on the modelization of Allee effect, fear effect, cannibalism and immigration. Some of the stated models show a very rich dynamics which is much more realistic that the one of the first historical models. 

The inclusion of Allee effect in a predator-prey system can have important effects on the dynamics; it can stabilize or destabilize the system or cause that the solutions of the system take a much longer time to reach a stable equilibrium point as seen in \cite{Allee1,Olivares2010}. 
We note that the Allee effect can produce different consequences depending on whether it is considered in the prey species or in the predator species, see  the works of Merdan \cite{Allee1} and Guan, Liu and Xie \cite{Allee2} and the comparative Figures \ref{fig:Allee1} and \ref{fig:Allee2}. Also limit cycles can appear in systems with Allee effect, as analyzed in the work of Wang, Xi and Wei \cite{Alleestrong}, and even two limit cycles can appear as proved in \cite{Olivares2011}.

The conclusions about the influence of fear effect are not the same in all the studied models. Considering a system with linear functional response and specialist predators, in \cite{fearWang} the authors show that fear does not affect the dynamics, but considering omnivorous predators, in \cite{fear} it is shown that the increase of the fear constant decreases the population density of the species in the singular point and can even cause extinction. Different behaviours for this two systems are compared in Figures \ref{fig:(25)(19)a} and \ref{fig:(25)(19)b}. Also for the system considered in \cite{fearWang}, it is proved that if we change the functional response from lineal to Holling type II, then the fear affects the dynamics: a value large enough for the fear parameter can help to maintain the asymptotic stability of a singular point.

The stabilizing effect of fear is also showed when combined with other biological characteristics. In \cite{miedocaza} the fear effect is combined with hunting cooperation and it is shown that when the system has oscillating behaviours, the increase of hunting cooperation or the increase of fear effect produces the stabilization. Also a stabilizing effect appears when combining the fear effect with prey refugee, as shown in \cite{fear-preyrefuge}.

In constrast, in \cite{Sasmal} we see that the stabilizing effect does not appear when fear is combined with Allee effect. In this case the only consequence is a reduction of the population density in the positive singular point. Similar conclusions are obtained in \cite{fear-alleeaditive} with a modified Allee effect.

Cannibalism had been principally considered in the predator species, and in that case the results agreed that cannibalism stabilizes the system and causes persistence in a population doomed to go extinct, but we have seen that considering cannibalism in prey species can produce the opposite result, as it can destabilize the system \cite{cannibalism2}. 

Regarding the presence of immigration we can observe that it has a stabilizing effect, reducing the amplitude of oscillations when they exist \cite{inm1}. Even very small immigration can cause that cyclic populations be stabilized, as proved in \cite{inm2}. Furthermore, this stabilization occurs whether immigration occurs on prey or predators.

We notice that the same methods and results are used in most papers for the study of the stability of the equilibria and the qualitative behaviour of the system. 
Principally, the authors use the standard linear stability analysis, they analyze the sign of the real part of the eigenvalues of the jacobian matrix of the system, sometimes by using the Routh-Hurwitz criterion. 

It is important to remember that for non-hyperbolic equilibria, the computation of the eigenvalues does not allow us to conclude on the stability. We have observed that in such cases the equilibria were not usually studied in the papers reviewed. Some other frequently used tools are the well-known Poincaré-Bendixon theorem, the Bendixon-Dulac criterion \cite{Libro} and the computation of Lyapunov coefficients \cite{Chicone}.

In most works, the authors make variable changes to obtain systems that are easier to study, for example by writing them down in a dimensionless form. In many cases, the systems are transformed into systems of Kolmogorov type, which are a generalization of Lotka-Volterra systems \eqref{LVgen} to an arbitrary degree. These systems are often used in biomathematics and their theoretical study  is currently progressing \cite{Kol1,Kol2,Kol3,Kol4,Kol5,Kol6,Kol7,Kol8}.
The transformation into Lienard type systems is also employed in some cases.

\vspace{0.5cm}
There are some open problems and some characteristics that would be interesting to investigate more in detail.
In the following we list some problems which researchers may find interesting to address in the future, which are proposed or motivated by the works included in this review.

\begin{itemize}
	\item Study new functions that can model the Allee effect, it is, functions representing the lowest population growth when population density is very low. This suggestion of considering different expression to model the same effect is also raised for some of the authors working on this topic, as for example in  \cite{Allee1}. We recall that according to biological facts, a function $\alpha(x)$ representing the Allee effect must verify that $\alpha'(x)$ is positive for any positive $x$ because the Allee effect decreases as density increases and $\lim_{x\to \infty}\alpha(x)=1$ because the Allee effect vanishes at high densities. \color{black}
	\vspace{0.15cm}
	
	\item   Among the systems that consider the Allee effect, the one proposed in \cite{Allee_genomn} is the one that shows a richer dynamics. As the authors say, it is possible that by changing the functional response and taking one of Holling type III, the behaviour is even more realistic, so it would be interesting to study it.
	\vspace{0.15cm}
	
	\item It would be interesting to consider models in which fear affects intraspecific competition, because as it is said in \cite{fearWang}, there are arguments in suport of this. It would be also important to get experimental evidence of this phenomenon to model it as realistic as possible.
	\vspace{0.15cm}
	
	\item In \cite{cannibalism2} it is proved that cannibalism can lead to limit cycle dynamics, but it would be interesting to determine the uniqueness or non uniqueness of limit cycles for cannibalistic populations, and answer the question of how many limit cycles can appear.
	\vspace{0.15cm}
	
	\item The effect of cannibalism on predator and prey simultaneously has not been studied
	as far as we know. But what happens when both are considered and predator cannibalism has a stabilizing effect while prey cannibalism has a destabilizing effect? A problem of this type is proposed in \cite{cannibalism2} for Holling Tanner models, but we think this question can be studied in any general model including cannibalism.
	\vspace{0.15cm}

	\item In \cite{inm2} system \eqref{41} was studied when there exists small immigration in prey or in predator, but the case with immigration in both species was not studied. It would be interesting to study the effect of immigration in both species simultaneously.
	\vspace{0.15cm}

	\item From a theoretical point of view, it would be specially interesting to carry out a study of the singular points of the systems for the values of the parameters for which they are non-hyperbolic, as in most articles this cases are omitted.
	\vspace{0.15cm}
	
	\item In \cite{Allee_genomn} the existence of Hopf bifurcations and Bogdanov-Takens bifurcations is proved numerically. It would be interesting to find an analytical proof of these results. 
	\vspace{0.15cm}
	
	\item  There are lots of effects and characteristics of the species and the environments that would be interesting to study, as for example, the consequences of anomalous diffusion and heterogeneous environments, see \cite{difusion}
	.

\end{itemize}

\vspace{6pt} 



\section{Contributions}
	All the authors have contributed in a similar way to this paper. All authors have read and agreed to the published version of the manuscript.

\section{Funding}The authors are partially supported by the Ministerio de Economía, Industria y Competitividad, Agencia Estatal de Investigación (Spain), grant MTM2016-79661-P (European FEDER support included, UE) and the Consellería de Educación, Universidade e Formación Profesional (Xunta de Galicia), grant ED431C 2019/10 with FEDER funds. The first author is also supported by the Ministerio de Educacion, Cultura y Deporte de España, contract FPU17/02125.









%

\end{document}